\documentclass[11pt]{article}
\usepackage{amsmath}
\usepackage{amsfonts}
\usepackage[mathscr]{euscript}
\usepackage{color}
\usepackage{graphics,graphicx}
\usepackage{subfig,float}
\usepackage{cancel}
\newtheorem{theo}{Theorem}[section]
\newtheorem{lem}[theo]{Lemma}

\newtheorem{defi}[theo]{Definition}

\newcommand{\mysection}[1]{\section{#1} \setcounter{equation}{0}}
\newcommand{\proof}{{\sc Proof.} \quad}
\newcommand{\proofc}{{\sc Proof} \ }
\newcommand{\be}{\begin{equation} \label}
\newcommand{\ee}{\end{equation}}
\newcommand{\bea}{\begin{eqnarray}\label}
\newcommand{\eea}{\end{eqnarray}}
\newcommand{\bas}{\begin{eqnarray*}}
\newcommand{\eas}{\end{eqnarray*}}
\newcommand{\bit}{\begin{itemize}}
\newcommand{\eit}{\end{itemize}}
\newcommand{\qed}{\hfill$\Box$ \vskip.2cm}
\newcommand{\nn}{\nonumber}
\newcommand{\R}{\mathbb{R}}
\newcommand{\N}{\mathbb{N}}
\newcommand{\pO}{\partial\Omega}
\newcommand{\bom}{\overline{\Omega}}
\newcommand{\eps}{\varepsilon}

\newcommand{\wsto}{\stackrel{\star}{\rightharpoonup}}
\newcommand{\hra}{\hookrightarrow}
\newcommand{\io}{\int_\Omega}

\newcommand{\abs}{\\[5pt]}

\newcommand{\ueps}{u_\eps}
\newcommand{\veps}{v_\eps}
\newcommand{\ts}{t_\star}
\begin{document}
\title{Relaxation by nonlinear diffusion enhancement in a\\ 
two-dimensional cross-diffusion model for urban crime propagation}
\author{
Nancy Rodriguez\footnote{rodrign@colorado.edu} \\
{\small Department of Applied Mathematics,}\\
{\small University of Colorado, Boulder, } \\
{\small 11 Engineering Dr., Boulder, Colorado 80309, USA}\\
\and
Michael Winkler\footnote{michael.winkler@math.uni-paderborn.de}\\
{\small Institut f\"ur Mathematik, Universit\"at Paderborn,}\\
{\small 33098 Paderborn, Germany} 
}

\date{}
\maketitle
\begin{abstract}
\noindent 
  We consider a class of macroscopic models for the spatio-temporal evolution of urban crime, as originally going back to
  Short et al. ({\it Math. Mod.~Meth.~Appl.~Sci.} {\bf 18}, 2008).
  The focus here is on the question how far a certain nonlinear enhancement in the random diffusion of criminal agents may exert visible
  relaxation effects.
  Specifically, in the context of the system
  \bas
   	\left\{ \begin{array}{l}
	u_t = \nabla \cdot (u^{m-1} \nabla u) - \chi \nabla \cdot \Big(\frac{u}{v} \nabla v \Big) - uv + B_1(x,t), \\[1mm]
	v_t = \Delta v +uv - v + B_2(x,t),
	\end{array} \right.
  \eas
  it is shown that whenever $\chi>0$ and the given nonnegative source terms $B_1$ and $B_2$ are sufficiently regular, the assumption
  \bas
	m>\frac{3}{2}
  \eas
  is sufficient to ensure that
  a corresponding Neumann-type initial-boundary value problem, posed in a smoothly bounded planar convex domain, admits locally bounded
 solutions for a wide class of arbitrary initial data. Furthermore, this solution is seen to be globally bounded if
  both $B_1$ and $B_2$ are bounded and $\liminf_{t\to\infty} \io B_2(\cdot,t)$ is positive.\abs
  This is supplemented by numerical evidence which, besides illustrating associated smoothing effects in particular situations
  of sharply structured initial data in the presence of such porous medium type diffusion mechanisms, 
  indicates a significant tendency toward support of singular structures in the linear diffusion case $m=1$.\abs
\end{abstract}
\newpage
\section{Introduction}\label{intro}
This manuscript is concerned with an adaptation a macroscopic model for the dynamics of urban crime, such as residential burglaries.
In its original version, as proposed by Short and collaborators in \cite{Short2008}, this model takes the form
\be{short}
   	\left\{ \begin{array}{ll}
	u_t = \nabla \cdot (D\nabla u) - 2 \nabla \cdot \Big(\frac{u}{v} \nabla v \Big) - uv + B_1(x,t),
	\qquad & x\in\Omega, \ t>0, \\[1mm]
	v_t = \Delta v +uv - v + B_2(x,t), 
	\qquad & x\in\Omega, \ t>0, \\[1mm]
    	\end{array} \right.
\ee
where
$u(x,t)$ denotes the density of criminals at location $x$ and time $t$ and $v(x,t)$ denotes the attractiveness field, which measures
the appeal of a location $x$ to a criminal agent at time $t$.   
The derivation of system \eqref{short} from first-principles heavily relies on assuming the routine activity theory \cite{Felson1987}, which asserts that crime revolves around three factors: a potential offender (criminal agents), 
a suitable target, and the absence of guardianship (e.g. police officers) (\cite{Felson1987}).  
This assumption thus provides a simple framework that allows one to  
move from an agent-based system to a continuum model.    
The other key assumption taken is the repeat and near-repeat victimization effect, stating that criminal activity 
in a certain location increases the probability of crime occurring again at the same or nearby locations.   
This effect has been observed in real-life data for 
crimes like residential burglaries (\cite{Johnson1997}, \cite{Short2009}).  
Essentially, this means that crime is self-exciting and in this framework it is modeled by the assumption that 
crime increases the attractiveness field.  In system \eqref{short} the parameter 
$D$ represents the effect that a crime that occurs in some location will have on immediate neighboring locations.\abs
The agent-based model introduced in \cite{Short2008} assumes that the probability that a crime occurs at a certain place and time follows a Poisson process with $v$ as its expectation.  
In the continuum limit given by \eqref{short} this implies that the number of crimes is proportional the attractiveness; thus, the expected number of crimes is given by the product $uv$.  
From the repeat victimization effect, we see that the attractiveness value increases with the number of crimes, giving rise to the 
summand $+uv$ in the $v$ equation.  Moreover, when a criminal agent commits a crime, 
it is supposed that they 
exit the system, giving rise to the term $-uv$ in the first equation.
This is based on the assumption that criminals want to keep a low-profile after committing a crime.  
Note that this hypothesis precludes us from considering serial crimes, where an agent might
rob various locations sequentially.    
To counteract the exit of criminal agents, the density of $u$ has a growth-pattern given by the known function $B_1(x,t).$ In \cite{Short2008} the growth-pattern
analyzed is given by a constant.  The near-repeat victimization effect is observed in the diffusion of the attractiveness value, the term $D\Delta v$ in the $v$ equation.  Finally, we note that in the context of (\ref{short}), 
criminals are assumed to move with a combination of unconditional dispersal, 
leading to the choice $D\equiv const.$ and hence, essentially, to the linear diffusion term $\Delta u$ in the first equation,
and conditional dispersal, $- 2 \nabla \cdot \Big(\frac{u}{v} \nabla v \Big)$, biased by high values of attractiveness.  \abs
The adaptation we consider in this work is based on the premise that criminals might have a tendency to avoid regions with a high density of other criminals.  
This seems reasonable in the cases when criminals want to avoid competition, or even suspect
that hotspot policing is being employed (\cite{Short2010}).  Hotspot policing is a strategy where the police 
force is deployed to areas with high crime (or equivalently a high density of criminals).  
In such cases, criminals are interested in avoiding areas with a high police density (or equivalently
areas with a high criminal density).  
The assumption that 
criminal agents tend to avoid police officers is a natural consequence of routine activity theory.  Indeed, one of the factors needed for crime to occur, based on this theory, is absence of guardianship.  Thus, 
criminal agents will not commit a crime in locations where there are police agents, and will instead choose to move away 
from areas with a high density of police. \abs

A natural approach to incorporate such a change in the movement strategy seems to consist in allowing the diffusion rate $D$ to depend
on $u$ and, in particular, to increase with $u$. Concentrating here on the apparently most prototypical algebraic choice of $D=D(u)$,
hence leading to porous medium type diffusion operators,
in the framework of a full no-flux initial-boundary value problem 
we subsequently consider the variant of (\ref{short}) given by
\be{0}
   	\left\{ \begin{array}{ll}
	u_t = \nabla \cdot (u^{m-1} \nabla u) - \chi \nabla \cdot \Big(\frac{u}{v} \nabla v \Big) - uv + B_1(x,t),
	\qquad & x\in\Omega, \ t>0, \\[1mm]
	v_t = \Delta v +uv - v + B_2(x,t), 
	\qquad & x\in\Omega, \ t>0, \\[1mm]
	\frac{\partial u}{\partial\nu}
	= \frac{\partial v}{\partial\nu}=0,
	\qquad & x\in\partial\Omega, \ t>0, \\[1mm]
	u(x,0)=u_0(x), \quad v(x,0)=v_0(x), 
	\qquad & x\in\Omega,
    	\end{array} \right.
\ee
in a bounded domain $\Omega\subset \R^2$ with smooth boundary. Here,
$B_1$ and $B_2$ are suitably regular nonnegative functions on $\Omega\times (0,\infty)$, $m>1$ is a given parameter and
$\chi$ is allowed to attain any positive value, thus including the choice $\chi=2$ in (\ref{short}) as a special case.\abs
We note that in order to keep the modeling framework as simple as possible,
in this work we do not independently model the dynamics of the police force by, e.g., describing their population density 
through an additional variable, but rather we make the simplifying assumption that the police force will match those of 
the criminal agents.

{\bf Main results: Blow-up suppression by strong diffusion enhancement.} \quad
Due to the potentially substantial destabilizing character of the self-enhanced cross-diffusive interaction therein,
systems of the form (\ref{short}) seem to bring about significant challenges
already at the level of basic solution theories. 
Accordingly, the few analytical findings available for (\ref{short}) and related systems are either restricted
to spatially one-dimensional settings (\cite{rodriguez_win1}), or address ranges of suitably small $\chi$ which do not contain
the relevant choice $\chi=2$ (\cite{freitag}), or concentrate 
on certain small-data solutions in cases of sufficiently small $B_1$ and $B_2$ (\cite{taowin_crime}), or resort to strongly 
generalized concepts of solvability which do not a priori preclude the emergence of singularities within finite time
(\cite{heihoff}, \cite{win_ANIHPC}). 
Although apparently no analytical study has rigorously detected the occurrence of such phenomena yet, the outcome of
numerical experiments supports the conjecture that indeed the linear diffusion mechanism in (\ref{short}) is insufficient to rule out
the possibility of explosions (cf.~also Section \ref{sect_numerics}).\abs
In contrast to this, we shall see that the presence of suitably strong nonlinear diffusion enhancement entirely suppresses
any such singular behavior in (\ref{0}) within finite time intervals, 
as expressed in the following statement on global existence of locally bounded solutions:
\begin{theo}\label{theo68}
  Let $\Omega\subset\R^2$ be a bounded convex domain with smooth boundary, and suppose that $\chi>0$, that
  \be{B}
	\mbox{$B_1$ and $B_2$ are nonnegative functions from $C^1(\bom\times [0,\infty))$,}
	\tag{B}
  \ee
  and that
  \be{m}
	m>\frac{3}{2}.
  \ee
  Then for any choice of functions $u_0$ and $v_0$ which are such that
  \be{init}
	\left\{ \begin{array}{l}
	u_0 \in W^{1,\infty}(\Omega) \quad \mbox{is nonnegative, and that} \\[1mm]
	v_0 \in W^{1,\infty}(\Omega) \quad \mbox{is positive in $\bom$,}
	\end{array} \right.
  \ee
  the problem (\ref{0}) possesses
  at least one global weak solution $(u,v)$ in the sense of Definition \ref{dw} below. 
  This solution is locally bounded in that
  \bas
	{\rm ess} \! \sup_{\hspace*{-4mm} t\in (0,T)} \|u(\cdot,t)\|_{L^\infty(\Omega)} <\infty
	\qquad \mbox{for all } T>0
  \eas
  and
  \bas
	{\rm ess} \! \sup_{\hspace*{-4mm} t\in (0,T)} \|v(\cdot,t)\|_{W^{1,q}(\Omega)} <\infty
	\qquad \mbox{for all $T>0$ and } q>2.
  \eas
\end{theo}
Under quite mild additional assumptions on $B_1$ and $B_2$, particularly fulfilled by any nonnegative 
$B_1=B_1(x)\in C^1(\bom)$ and $0\not\equiv B_2=B_2(x)\in C^1(\bom)$, solutions can be found which are in fact globally bounded,
meaning that in such cases moreover even any infinite-time singularity formation is ruled out:
\begin{theo}\label{theo69}
  Assume that $\Omega\subset\R^2$ be a bounded convex domain with smooth boundary, that $m>\frac{3}{2}$ and $\chi>0$,
  and that $(u_0,v_0)$ satisfies (\ref{init}), and suppose furthermore that $B_1$ and $B_2$ are such that beyond
  (\ref{B}) we have
  \be{B1}
	\sup_{(x,t) \in \Omega\times (0,\infty)} \Big\{ B_1(x,t) + B_2(x,t)\Big\} <\infty
	\tag{B1}
  \ee
  and
  \be{B2}
	\liminf_{t\to\infty} \io B_2(x,t) dx >0.
	\tag{B2}
  \ee
  Then (\ref{0}) admits a global weak solution according to Definition \ref{dw} which is globally bounded in the sense 
  that
  \be{69.1} 
	{\rm ess} \! \sup_{\hspace*{-4mm} t>0} \|u(\cdot,t)\|_{L^\infty(\Omega)} <\infty
  \ee
  and
  \be{69.2}
	{\rm ess} \! \sup_{\hspace*{-4mm} t>0} \|v(\cdot,t)\|_{W^{1,q}(\Omega)} <\infty
	\qquad \mbox{for all } q>2.
  \ee
\end{theo}
Accompanied and illustrated by outcomes of corresponding numerical simulations, to be presented in Section \ref{sect_numerics},
these results quantitatively identify an effect of the considered diffusion strengthening on overcrowding prevention.
This seems to indicate that nonlinear migration mechanisms of the said flavor	
may stabilize systems of the considered form by precluding a model breakdown due to the emergence of singularities.  
Viewed in the contexts of the addressed application seems to be of relevance, especially due 
to the nontrivial size of criminal agents.
As partially seen in Section \ref{sect_numerics},
the description of crime hotspot formation, as known to occur in associated typical real-life situations, is thereby
transported to mathematical sceneries involving structured but bounded spatial profiles, rather than exploding solutions
such as naturally going along with Keller-Segel type modeling of aggregation in populations of microbial individuals
(\cite{herrero_velazquez}, \cite{win_JMPA}; see also \cite{BBTW}). 
\mysection{Regularization and basic properties}
In order to conveniently regularize (\ref{0}), we combine the essence of the corresponding procedure in \cite{win_ANIHPC} with
a standard non-degenerate approximation of porous medium type diffusion operators, and hence we shall subsequently consider the problems
\be{0eps}
   	\left\{ \begin{array}{ll}
	u_{\eps t} = \nabla \cdot \big( (\ueps+\eps)^{m-1} \nabla\ueps\big)
	- \chi \nabla\cdot \Big( \frac{\ueps}{\veps} \nabla\veps\Big) - \ueps\veps + B_1(x,t),
	\qquad & x\in\Omega, \ t>0, \\[1mm]
	v_{\eps t} = \Delta \veps + \frac{\ueps\veps}{1+\eps\ueps\veps} - \veps + B_2(x,t), 
	\qquad & x\in\Omega, \ t>0, \\[1mm]
	\frac{\partial \ueps}{\partial\nu}
	= \frac{\partial \veps}{\partial\nu}=0,
	\qquad & x\in\partial\Omega, \ t>0, \\[1mm]
	\ueps(x,0)=u_0(x), \quad \veps(x,0)=v_0(x), 
	\qquad & x\in\Omega,
    	\end{array} \right.
\ee
for $\eps\in (0,1)$, which indeed are all globally solvable in the classical sense:
\begin{lem}\label{lem_loc}
  Assume (\ref{B}) and (\ref{init}), and let $m>1$ and $\eps\in (0,1)$. Then there exist functions
  \bas
	\left\{ \begin{array}{l}
	\ueps \in C^0(\bom\times [0,\infty)) \cap C^{2,1}(\bom\times (0,\infty)), \\[1mm]
	\veps \in \bigcap_{p>2} C^0([0,\infty); W^{1,p}(\Omega)) \cap C^{2,1}(\bom\times (0,\infty)),
	\end{array} \right.
  \eas
  which solve (\ref{0eps}) classically in $\bom\times [0,\infty)$, and which are such that $\ueps>0$ 
  in $\bom\times (0,\infty)$ and $\veps>0$ in $\bom\times [0,\infty)$.
\end{lem}
\proof
  This can be seen by a straightforward adaptation of the reasoning in \cite[Section 2]{win_ANIHPC} on the basis
  of standard results on local existence and extensibility, as provided e.g.~by the general theory in \cite{amann}.
\qed
Throughout the sequel, without further explicit mentioning we shall assume that (\ref{B}) and (\ref{init}) are satisfied, 
and for $m>1$ and $\eps\in (0,1)$ we let $(\ueps,\veps)$ denote the solutions of (\ref{0eps}) gained above.\abs
In our respective formulation of statements on regularity of these solutions, we find it convenient to make use of the following
notational convention concerning a certain time independence of constants under the hypotheses (\ref{B1}) and (\ref{B2}).
\begin{defi}\label{dk}
  Let $K: (0,\infty) \to (0,\infty)$. We then say that $K$ satisfies (K) if $K$ has the property that
  \bas
	\sup_{T>0} K(T)<\infty
	\qquad \mbox{whenever (\ref{B1}) and (\ref{B2}) hold.}
  \eas
\end{defi}
With reference to this property, our first basic statement on a pointwise lower bound for the second solution component,
resembling similar information found in \cite{rodriguez_win1} and \cite{win_ANIHPC} already, reads as follows.
\begin{lem}\label{lem01}
  Let $m>1$. Then
  there exists $K: (0,\infty) \to (0,\infty)$ fulfilling (K) such that whenever $T>0$,
  \be{01.1}
	\frac{1}{\veps(x,t)} \le K(T)
	\qquad \mbox{for all $x\in\Omega, t\in (0,T)$ and } \eps\in (0,1).
  \ee
\end{lem}
\proof
  Firstly, in view of the nonnegativity of $\ueps, \veps$ and $B_2$ it follows by a comparison argument that
  \be{01.2}
	\veps(x,t) \ge \bigg\{ \inf_{y\in\Omega} v_0(y) \bigg\} \cdot e^{-t}
	\qquad \mbox{for all $x\in\Omega, t>0$ and } \eps\in (0,1).
  \ee
  Moreover, the convexity of $\Omega$ allows us to import from \cite{fujie_diss} 
  a result on a pointwise positivity feature of the Neumann heat semigroup
  $(e^{t\Delta})_{t\ge 0}$ on $\Omega$ to fix $c_1>0$ fulfilling
  \bas
	e^{t\Delta} \psi \ge c_1 \io \psi
	\quad \mbox{in $\Omega$} \qquad
	\mbox{for all $t>1$ and any nonnegative } \psi\in C^0(\bom),
  \eas
  whence again by the comparison principle, for arbitrary $t_0\ge 0$ we can estimate
  \bea{01.3}
	\veps(\cdot,t)
	&=& e^{t(\Delta-1)} v_0
	+ \int_0^t e^{(t-s)(\Delta-1)} \Big\{ \frac{\ueps(\cdot,s)\veps(\cdot,s)}{1+\eps\ueps(\cdot,s)\veps(\cdot,s)} + B_2(\cdot,s) 
		\Big\} ds \nn\\
	&\ge& c_1 \int_{t_0}^t e^{-(t-s)} \cdot \bigg\{ \io B_2(\cdot,s) \bigg\} ds \nn\\
	&\ge& c_1 \cdot \bigg\{ \inf_{s>t_0} \io B_2(\cdot,s)\bigg\} \cdot \int_{t_0}^t e^{-(t-s)} ds \nn\\
	&=& c_1 \cdot \bigg\{ \inf_{s>t_0} \io B_2(\cdot,s)\bigg\} \cdot (1-e^{-(t-t_0)}) \nn\\
	&\ge& (1-e^{-1}) c_1 \cdot \bigg\{ \inf_{s>t_0} \io B_2(\cdot,s)\bigg\}
	\quad \mbox{in } \Omega 
	\qquad \mbox{for all } t>t_0+1.
  \eea
  Combining (\ref{01.2}) with (\ref{01.3}) readily yields (\ref{01.1}) with some $K$ satisfying (K).
\qed
Likewise, our second basic observation has quite closely related precedents in
\cite{rodriguez_win1} and \cite{win_ANIHPC}.
\begin{lem}\label{lem1}
  Let $m>1$. Then
  there exists $K: (0,\infty) \to (0,\infty)$ such that (K) holds, and such that for all $T>0$,
  \be{1.1}
	\io \ueps(\cdot,t) \le K(T)
	\qquad \mbox{for all $t\in (0,T)$ and } \eps\in (0,1),
  \ee
  and that
  \be{1.2}
	\io \veps(\cdot,t) \le K(T)
	\qquad \mbox{for all $t\in (0,T)$ and } \eps\in (0,1).
  \ee
\end{lem}
\proof
  According to Lemma \ref{lem01}, we can find $k_1: (0,\infty)\to (0,\infty)$ with the corresponding property (K) such that for
  all $T>0$,
  \be{1.3}
	\frac{1}{\veps} \le k_1(T)
	\quad \mbox{in } \Omega\times (0,T)
	\qquad \mbox{for all } \eps\in (0,1).
  \ee
  Then letting 
  \be{1.33}
	k_2(T):=\min \Big\{ 1 \, , \, \frac{1}{2k_1(T)} \Big\},
	\qquad T>0,
  \ee
  we use (\ref{0eps}) to see that given any $T>0$, for all $t>0$ and each $\eps\in (0,1)$ we have
  \bea{1.4}
	\hspace*{-20mm}
	\frac{d}{dt} \bigg\{ 2\io \ueps + \io \veps \bigg\}
	&+& k_2(T) \cdot \bigg\{ 2\io \ueps + \io \veps \bigg\} \nn\\
	&=& - 2 \io \ueps\veps + 2 \io B_1 \nn\\
	& & - \io \veps + \io \frac{\ueps\veps}{1+\eps\ueps\veps} + \io B_2 \nn\\
	& & + 2k_2(T) \io \ueps + k_2(T) \io \veps \nn\\
	&\le& - \io \ueps\veps 
	+ 2k_2(T) \io \ueps
	+ 2\io B_1 + \io B_2,
  \eea
  because $k_2(T) \le 1$.
  Using that moreover $2k_2(T) \le \frac{1}{k_1(T)}$ and hence
  \bas
	- \io \ueps\veps + 2k_2(T) \io \ueps
	\le - \frac{1}{k_1(T)} \io \ueps + 2k_2(T) \io \ueps \le 0
	\qquad \mbox{for all $t\in (0,T)$ and } \eps\in (0,1)
  \eas
  by (\ref{1.3}), from (\ref{1.4}) we infer that
  \bas
	& & \hspace*{-20mm}
	\frac{d}{dt} \bigg\{ 2\io \ueps + \io \veps \bigg\}
	+ k_2(T) \cdot \bigg\{ 2\io \ueps + \io \veps \bigg\} \nn\\
	&\le& k_3(T):=\sup_{s\in (0,T)} \bigg\{ 2\io B_1 + \io B_2 \bigg\}
	\qquad \mbox{for all $t\in (0,T)$ and } \eps\in (0,1).
  \eas
  Therefore, an ODE comparison shows that
  \bas
	2 \io \ueps(\cdot,t) + \io \veps(\cdot,t)
	\le \max \bigg\{ 2\io u_0 + \io v_0 \, , \, \frac{k_3(T)}{k_2(T)} \bigg\}
	\qquad \mbox{for all $t\in (0,T)$ and } \eps\in (0,1),
  \eas
  from which both (\ref{1.1}) and (\ref{1.2}) result upon the observation that in view of (\ref{1.3}), (K) holds for the
  function $\frac{k_3}{k_2}$.
\qed
\mysection{Estimates for $\veps$ in $W^{1,q}(\Omega)$ with $q\le 2$}
The following estimate essentially reproduces a similar finding from \cite[Lemma 3.1]{win_ANIHPC} to the present framework
involving slightly different hypotheses on $B_1$ and $B_2$.
\begin{lem}\label{lem2}
  Assume that $m>1$, and let $p\in (0,1)$. 
  Then there exists a function $K\equiv K^{(p)}: (0,\infty) \to (0,\infty)$ which satisfies (K) and is such that
  whenever $T>0$,
  \be{2.1}
	\int_t^{t+1} \io \veps^{p-2} |\nabla\veps|^2 \le K(T)
	\qquad \mbox{for all $t\in (0,T)$ and } \eps\in (0,1).	
  \ee
\end{lem}
\proof
  Relying on Lemma \ref{lem1}, we can fix a mapping $k_1: (0,\infty)\to (0,\infty)$ which enjoys the boundedness feature in (K)
  and is such that for all $T>0$,
  \bas
	\io \veps(\cdot,t) \le k_1(T)
	\qquad \mbox{for all $t\in (0,T+1)$ and } \eps\in (0,1),
  \eas
  whence given $p\in (0,1)$ we can use Young's inequality to see that
  \be{2.2}
	\io \veps^p(\cdot,t) 
	\le \io \Big(\veps(\cdot,t)+1\Big)
	\le k_1(T) +|\Omega|
	\qquad \mbox{for all $t\in (0,T+1)$ and } \eps\in (0,1).
  \ee
  Since according to (\ref{0eps}) we have
  \bas
	\frac{1}{p} \frac{d}{dt} \io \veps^p
	&=& (1-p) \io \veps^{p-2} |\nabla\veps|^2
	- \io \veps^p + \io \frac{\ueps \veps^p}{1+\eps\ueps} + \io \veps^{p-1} B_2 \\
	&\ge& (1-p) \io \veps^{p-2} |\nabla\veps|^2
	- \io \veps^p 
	\qquad \mbox{for all $t>0$ and } \eps\in (0,1)
  \eas
  and thus
  \bas
	(1-p)\int_t^{t+1} \io \veps^{p-2} |\nabla\veps|^2
	&\le& \frac{1}{p} \io \veps^p(\cdot,t+1) 
	- \frac{1}{p} \io \veps^p(\cdot,t)
	+ \int_t^{t+1} \io \veps^p \\
	&\le& \frac{1}{p} \io \veps^p(\cdot,t+1)
	+ \int_t^{t+1} \io \veps^p
	\qquad \mbox{for all $t>0$ and } \eps\in (0,1),
  \eas
  utilizing (\ref{2.2}) to estimate
  \bas
	\frac{1}{p} \io \veps^p(\cdot,t+1)
	+ \int_t^{t+1} \io \veps^p
	\le \frac{1}{p} \cdot \big( k_1(T)+|\Omega|\big) + k_1(T) + |\Omega|
	\qquad \mbox{for all $t\in (0,T)$ and } \eps\in (0,1)
  \eas
  we readily arrive at (\ref{2.1}) upon an evident choice of $K$.
\qed
Besides being of independent use in some of our subsequent estimates (see Lemma \ref{lem43} and Lemma \ref{lem52}),
Lemma \ref{lem2}, through suitable interpolation involving Lemma \ref{lem1}, also entails the following boundedness property
of $\veps$ with respect to the norm in $W^{1,q}(\Omega)$ for $q\in [1,2)$ arbitrarily close to $2$.
\begin{lem}\label{lem3}
 Suppose that $m>1$ and let $q\in [1,2)$. 
  Then there exists $K\equiv K^{(q)}: (0,\infty) \to (0,\infty)$ fulfilling (K) with the property that
  for all $T>0$, any $\eps\in (0,1)$ and each $t\in (0,T)$ fulfilling $t\ge 2$ one can find $t_0=t_0(t,\eps) \in (t-2,t-1)$ such that
  \be{3.1}
	\|\veps(\cdot,t_0)\|_{W^{1,q}(\Omega)} \le K(T).
  \ee
\end{lem}
\proof
  We evidently need to define $K(T)$ for $T\ge 2$ only, and to achieve this we first employ Lemma \ref{lem2} and Lemma \ref{lem1}
  to find $k_i: (0,\infty)\to (0,\infty)$, $i\in \{1,2\}$, which comply with (K) and are such that whenever $T\ge 2$,
  \be{3.2}
	\int_{t-2}^{t-1} \io \veps^{-\frac{3}{2}}|\nabla\veps|^2 \le k_1(T)
	\qquad \mbox{for all $t \in [2,T]$ and } \eps\in (0,1)
  \ee
  and
  \be{3.3}
	\io \veps \le k_2(T)
	\qquad \mbox{for all $t \in (0,T)$ and } \eps\in (0,1).
  \ee
  Moreover, given $q\in [1,2)$ we define $p=p(q):=\frac{3q}{2(2-q)}>1$ and make use of the continuity of the embedding
  $W^{1,2}(\Omega) \hra L^{4p}(\Omega)$ to fix $c_1=c_1(q)>0$ such that
  \be{3.4}
	\|\varphi\|_{L^{4p}(\Omega)}^{4p} \le c_1 \|\nabla\varphi\|_{L^2(\Omega)}^{4p} 
	+ c_1\|\varphi\|_{L^4(\Omega)}^{4p}
	\qquad \mbox{for all } \varphi\in W^{1,2}(\Omega).
  \ee
  Now letting $T\ge 2$ and $t\in [2,T]$ be arbitrary, from (\ref{3.2}) we infer the existence of $t_0=t_0(t,\eps)\in (t-2,t-1)$ 
  such that
  \be{3.5}
	\io \veps^{-\frac{3}{2}}(\cdot,t_0) |\nabla\veps(\cdot,t_0)|^2 \le k_1(T),
  \ee
  which in conjunction with (\ref{3.4}) and (\ref{3.3}) entails that
  \bas
	\io \veps^p(\cdot,t_0)
	&=& \|\veps^\frac{1}{4}(\cdot,t_0)\|_{L^{4p}(\Omega)}^{4p} \\
	&\le& c_1 \|\nabla\veps^\frac{1}{4}(\cdot,t_0)\|_{L^2(\Omega)}^{4p}
	+ c_1\|\veps^\frac{1}{4}(\cdot,t_0)\|_{L^4(\Omega)}^{4p} \\
	&=& \frac{c_1}{4^{4p}} \cdot \bigg\{ \io \veps^{-\frac{3}{2}}(\cdot,t_0) |\nabla\veps(\cdot,t_0)|^2 \bigg\}^{2p}
	+ c_1 \cdot \bigg\{ \io \veps(\cdot,t_0) \bigg\}^p \\
	&\le& k_3(T)\equiv k_3^{(q)}(T)
	:=\frac{c_1}{4^{4p}} \cdot k_1^{2p}(T) + c_1 k_p^2(T).
  \eas
  Once more combined with (\ref{3.5}), due to Young's inequality and thanks to our definition of $p$ this shows that
  \bas
	\io |\nabla\veps(\cdot,t_0)|^q
	&=& \io \Big\{ \veps^{-\frac{3}{2}}(\cdot,t_0) |\nabla\veps(\cdot,t_0)|^2 \Big\}^\frac{q}{2} 
		\cdot \veps^\frac{3q}{4}(\cdot,t_0) \\
	&\le& \io \veps^{-\frac{3}{2}}(\cdot,t_0)|\nabla\veps(\cdot,t_0)|^2
	+ \io \veps^\frac{3q}{2(2-q)}(\cdot,t_0) \\[1mm]
	&\le& k_1(T) + k_3(T).
  \eas
  In view of (\ref{3.3}), this implies the claimed boundedness property in $W^{1,q}(\Omega)$.
\qed
\mysection{Superlinear integrability properties of $\ueps$}
Our derivation of further regularity properties of $\veps$ will crucially rely on the following a priori information 
on the first solution component, obtained by means of a standard testing procedure on the basis of Lemma \ref{lem2}.
\begin{lem}\label{lem43}
  Let $m>1$. Then there exists $K: (0,\infty) \to (0,\infty)$ satisfying (K) such that if $T>0$ then
  \be{43.2}
	\left\{ \begin{array}{ll}
	\displaystyle
	\int_t^{t+1} \io (\ueps+\eps)^{m-1} (\ueps+1)^{m-3} |\nabla\ueps|^2 \le K(T)
	\quad \mbox{for all $t\in (0,T)$ and $\eps\in (0,1)$}
	\qquad & \mbox{if } m\in (1,2], \\[4mm]
	\displaystyle
	\int_t^{t+1} \io (\ueps+\eps)^{2m-4} |\nabla\ueps|^2 \le K(T)
	\quad \mbox{for all $t\in (0,T)$ and $\eps\in (0,1)$}
	\qquad & \mbox{if } m>2.
	\end{array} \right.
  \ee
\end{lem}
\proof
  We again return to Lemma \ref{lem01} and additionally employ Lemma \ref{lem2} to pick $k_i: (0,\infty) \to (0,\infty)$, $i\in \{1,2\}$,
  which comply with (K) and are such that if $T>0$ then
  \be{43.4}
	\veps \ge \frac{1}{k_1(T)} 
	\quad \mbox{in } \Omega\times (0,T)
	\qquad \mbox{for all } \eps\in (0,1)
  \ee
  and
  \be{43.44}
	\int_t^{t+1} \io \veps^{-\frac{3}{2}}|\nabla\veps|^2 \le k_2(T)
	\qquad \mbox{for all $t\in (0,T)$ and } \eps\in (0,1),
  \ee
  and to further prepare our argument for large $m$ we utilize Young's inequality together with (\ref{B1}) to find 
  $k_3: (0,\infty)\to (0,\infty)$ which is such that (K) holds and that if $m>2$, then whenever $T>0$,
  \be{43.5}
	\bigg\{ \frac{1}{k_1(T)} + \|B_1(\cdot,t)\|_{L^\infty(\Omega)} \bigg\} \cdot \xi^{m-2}
	\le \frac{1}{2k_1(T)} \xi^{m-1} + k_3(T)
	\qquad \mbox{for all } t\in (0,T).
  \ee
  Now for such $m$, we use $(\ueps+\eps)^{m-2}$ as a test function in (\ref{0eps}) to see that 
  \bea{43.55}
	& & \hspace*{-20mm}
	\frac{1}{m-1} \frac{d}{dt} \io (\ueps+\eps)^{m-1}
	+ (m-2) \io (\ueps+\eps)^{2m-4} |\nabla\ueps|^2 \nn\\
	&=& (m-2)\chi \io \ueps (\ueps+\eps)^{m-3} \nabla\ueps \cdot \frac{\nabla\veps}{\veps} \nn\\
	& & - \io \ueps (\ueps+\eps)^{m-2} \veps
	+ \io (\ueps+\eps)^{m-2} B_1
	\qquad \mbox{for all $t>0$ and } \eps\in (0,1),
  \eea
  where once more by Young's inequality, and by (\ref{43.4}), given $T>0$ we can estimate
  \bas
	& & \hspace*{-12mm}
	(m-2)\chi \io \ueps(\ueps+\eps)^{m-3} \nabla\ueps \cdot \frac{\nabla\veps}{\veps} \\
	&\le& \frac{m-2}{2} \io (\ueps+\eps)^{2m-4} |\nabla\ueps|^2 
	+ \frac{(m-2)\chi^2}{2} \io \Big(\frac{\ueps}{\ueps+\eps}\Big)^2 \frac{|\nabla\veps|^2}{\veps^2} \\
	&\le& \frac{m-2}{2} \io (\ueps+\eps)^{2m-4} |\nabla\ueps|^2 
	+ \frac{(m-2)\chi^2}{2} k_1^\frac{1}{2}(T) \io \veps^{-\frac{3}{2}}|\nabla\veps|^2
	\qquad \mbox{for all $t\in (0,T)$ and } \eps\in (0,1),
  \eas
  and where (\ref{43.4}) together with (\ref{43.5}) ensures that for all $T>0$,
  \bas
	& & \hspace*{-26mm}
	- \io \ueps(\ueps+\eps)^{m-2}\veps
	+ \io (\ueps+\eps)^{m-2} B_1 \\
	&\le& - \frac{1}{k_1(T)} \io \ueps (\ueps+\eps)^{m-2}
	+ \io (\ueps+\eps)^{m-2} B_1 \\
	&=& - \frac{1}{k_1(T)} \io (\ueps+\eps)^{m-1}
	+ \frac{\eps}{k_1(T)} \io (\ueps+\eps)^{m-2}
	+ \io (\ueps+\eps)^{m-2} B_1 \\
	&\le& - \frac{1}{k_1(T)} \io (\ueps+\eps)^{m-1}
	+ \bigg\{ \frac{1}{k_1(T)} + \|B_1(\cdot,t)\|_{L^\infty(\Omega)} \bigg\} \cdot \io (\ueps+\eps)^{m-2} \\
	&\le& - \frac{1}{2 k_1(T)} \io (\ueps+\eps)^{m-1}
	+ k_3(T) |\Omega|
	\qquad \mbox{for all $t\in (0,T)$ and } \eps\in (0,1).
  \eas
  Therefore, (\ref{43.55}) entails that whenever $T>0$,
  \bea{43.56}
	& & \hspace*{-20mm}
	\frac{1}{m-1} \frac{d}{dt} \io (\ueps+\eps)^{m-1}
	+ \frac{m-2}{2} \io (\ueps+\eps)^{2m-4} |\nabla\ueps|^2
	+ \frac{1}{2k_1(T)} \io (\ueps+\eps)^{m-1} \nn\\
	&\le& \frac{(m-2)\chi^2}{2} k_1^\frac{1}{2}(T) \io \veps^{-\frac{3}{2}}|\nabla\veps|^2
	+ k_3(T)|\Omega|
	\qquad \mbox{for all $t\in (0,T)$ and } \eps\in (0,1),
  \eea
  from which in light of (\ref{43.44}) the respective inequality in (\ref{43.2}) readily results upon an integration in time.\abs
  We are thus left with the case when $m\in (1,2]$, in which we let 
  $\Phi(\xi):=\int_0^{\xi} \int_0^{\sigma} (\tau+1)^{m-3}d\tau d\sigma$, $\xi\ge 0$, and noting that 
  $\Phi''(\xi)= (\xi+1)^{m-3}$ for all $\xi\ge 0$ we once more resort to (\ref{0eps})
  to see that similarly to the above, given an arbitrary $T>0$ we have
  \bea{43.6}
	\frac{d}{dt} \io \Phi(\ueps)
	&=& \io \Phi'(\ueps) \cdot 
	\bigg\{ \nabla \cdot \big( (\ueps+\eps)^{m-1} \nabla\ueps\big) - \chi \nabla\cdot \Big(\frac{\ueps}{\veps}\nabla\veps\Big)
	- \ueps\veps + B_1 \bigg\} \nn\\
	&=& - \io (\ueps+\eps)^{m-1} \Phi''(\ueps) |\nabla\ueps|^2
	+ \chi \io \ueps \Phi''(\ueps) \nabla\ueps \cdot \frac{\nabla\veps}{\veps} \nn\\
	& & - \io \ueps \Phi'(\ueps)\veps
	+ \io \Phi'(\ueps) B_1 \nn\\
	&=& - \io (\ueps+\eps)^{m-1} (\ueps+1)^{m-3} |\nabla\ueps|^2
	+ \chi \io \ueps(\ueps+1)^{m-3} \nabla\ueps \cdot \frac{\nabla\veps}{\veps} \nn\\
	& & - \io \ueps \Phi'(\ueps)\veps
	+ \io \Phi'(\ueps) B_1 \nn\\
	&\le& - \frac{1}{2} \io (\ueps+\eps)^{m-1} (\ueps+1)^{m-3} |\nabla\ueps|^2
	+ \frac{\chi^2}{2} \io \ueps^2 (\ueps+\eps)^{1-m} (\ueps+1)^{m-3} \frac{|\nabla\veps|^2}{\veps^2} \nn\\
	& & - \io \ueps \Phi'(\ueps)\veps
	+ \io \Phi'(\ueps) B_1 \nn\\
	&\le& - \frac{1}{2} \io (\ueps+\eps)^{m-1} (\ueps+1)^{m-3} |\nabla\ueps|^2
	+ \frac{\chi^2}{2} k_1^\frac{1}{2}(T) \io veps^{-\frac{3}{2}}|\nabla\veps|^2 \nn\\
	& & - \io \ueps \Phi'(\ueps)\veps
	+ \io \Phi'(\ueps) B_1 
	\qquad \mbox{for all $t\in (0,T)$ and } \eps\in (0,1)
  \eea
  because of the pointwise inequality 
  \bas
	\ueps^2 (\ueps+\eps)^{1-m} (\ueps+1)^{m-3} \le (\ueps+\eps)^{3-m} (\ueps+1)^{m-3} \le 1,
  \eas
  valid throughout $\Omega\times (0,\infty)$ for each $\eps\in (0,1)$ and any such $m$.
  Now from the definition of $\Phi$ we furthermore see that for all $\xi\ge 0$,
  \be{43.7}
	\Phi'(\xi) = \left\{ \begin{array}{ll}
	- \frac{1}{2-m} (\xi+1)^{m-2} + \frac{1}{2-m},
	\qquad & \mbox{if }, m\in (1,2), \\[2mm]
	\ln (\xi+1),
	\qquad & \mbox{if } m=2,
	\end{array} \right.
  \ee
  and
  \bas
	\Phi(\xi) = \left\{ \begin{array}{ll}
	- \frac{1}{(2-m)(m-1)} (\xi+1)^{m-1} + \frac{1}{2-m} \xi + \frac{1}{(2-m)(m-1)},
	\qquad & \mbox{if } m\in (1,2), \\[2mm]
	(\xi+1)\ln (\xi+1) -\xi,
	\qquad & \mbox{if } m=2,
	\end{array} \right.
  \eas
  from which it follows that for each $\xi\ge 0$,
  \bas
	\xi\Phi'(\xi) - \Phi(\xi) = \left\{ \begin{array}{ll}
	\frac{1}{m-1} (\xi+1)^{m-1} - \frac{1}{2-m} (\xi+1)^{m-2} - \frac{1}{(2-m)(m-1)},
	\qquad & \mbox{if } m\in (1,2), \\[2mm]
	\xi - (\xi+1)\ln (\xi+1),
	\qquad & \mbox{if } m=2.
	\end{array} \right.
  \eas
  Using that $(\xi+1)^{m-2} \le 1$ for all $\xi\ge 0$ when $m<2$, and that $\ln (\xi+1) \le \xi$ for any $\xi\ge 0$, we thus obtain
  $c_1>0$ and $c_2>0$ such that in both cases,
  \bas
	\xi\Phi'(\xi) - \Phi(\xi) \ge -c_1
	\qquad \mbox{for all } \xi\ge 0
  \eas
  and
  \bas
	\Phi'(\xi) \le \xi+c_2
	\qquad \mbox{for all } \xi\ge 0.
  \eas
  As (\ref{43.7}) moreover implies that $\Phi'\ge 0$ on $[0,\infty)$, once again going back to (\ref{43.4}) and recalling
  Lemma \ref{lem1}, we can thus find $k_4: (0,\infty) \to (0,\infty)$ fulfilling (K) such that
  for fixed $T>0$ we can estimate the two rightmost summands in (\ref{43.6}) according to
  \bas
	- \io \ueps \Phi'(\ueps)\veps
	+ \io \Phi'(\ueps) B_1 
	&\le& - \frac{1}{k_1(T)} \io \ueps \Phi'(\ueps)
	+ \|B_1(\cdot,t)\|_{L^\infty(\Omega)} \io \Phi'(\ueps) \\
	&\le& - \frac{1}{k_1(T)} \io \big\{ \Phi(\ueps)+c_1\big\}
	+ \|B_1(\cdot,t)\|_{L^\infty(\Omega)} \io (\ueps+c_2) \\
	&\le& - \frac{1}{k_1(T)} \io \Phi(\ueps) + k_4(T)
	\qquad \mbox{for all $t\in (0,T)$ and } \eps\in (0,1).
  \eas
  For any such $T$, from (\ref{43.6}) we consequently derive the analogue of (\ref{43.56}) given by
  \bas
	& & \hspace*{-30mm}
	\frac{d}{dt} \io \Phi(\ueps)
	+ \frac{1}{2} \io (\ueps+\eps)^{m-1} (\ueps+1)^{m-3} |\nabla\ueps|^2 
	+ \frac{1}{k_1(T)} \io \Phi(\ueps) \\
	&\le& \frac{\chi^2}{2} k_1^\frac{1}{2}(T) \io \veps^{-\frac{3}{2}}|\nabla\veps|^2
	+ k_4(T)
	\qquad \mbox{for all $t\in (0,T)$ and } \eps\in (0,1),
  \eas
  which due to (\ref{43.44}) and the evident nonnegativity of $\Phi$ entails the claimed inequality in (\ref{43.2}) also for such
  values of $m$.
\qed
An interpolation of the latter with the $L^1$ bound provided by Lemma \ref{lem1}, namely, yields a spatio-temporal
integral estimate for $\ueps$ itself which involves superlinear summability powers conveniently increasing with $m$.
\begin{lem}\label{lem44}
  Let $m>1$. Then
  there exists $K: (0,\infty) \to (0,\infty)$ satisfying (K) such that whenever $T>0$,
  \be{44.1}
	\int_t^{t+1} \io \ueps^{2m-1} \le K(T)
	\qquad \mbox{for all $t\in (0,T)$ and } \eps\in (0,1).
  \ee
\end{lem}
\proof
  We fix $\rho\in C^0([0,\infty))$ such that $\rho\equiv 0$ on $[0,1]$, $\rho(\xi)=\xi^{m-2}$ for all $\xi\ge 2$ and
  $0\le \rho(\xi) \le \xi^{m-2}$ for all $\xi\ge 0$, and let $P(\xi):=\int_0^\xi \rho(\sigma)d\sigma$ for $\xi\ge 0$.
  Then $P$ belongs to $C^1([0,\infty))$ and satisfies $P(\xi) \le \frac{\xi^{m-1}}{m-1}$ as well as
  \be{44.3}
	P(\xi) \ge \int_2^\xi \rho(\sigma) d\sigma
	= \frac{\xi^{m-1}-2^{m-1}}{m-1}
	\ge c_1 \xi^{m-1}
	\qquad \mbox{for all } \xi\ge 3
  \ee
  with $c_1:=\frac{1-(\frac{2}{3})^{m-1}}{m-1}>0$.
  Since thus
  \bas
	\|P(\ueps)\|_{L^\frac{1}{m-1}(\Omega)}^\frac{1}{m-1}
	\le (m-1)^{-\frac{1}{m-1}} \io \ueps
	\qquad \mbox{for all $t>0$ and } \eps\in (0,1)
  \eas
  and
  \bas
	\|\nabla P(\ueps)\|_{L^2(\Omega)}^2
	= \io \rho^2(\ueps)|\nabla\ueps|^2
	\le \int_{\{\ueps\ge 1\}} \ueps^{2m-4} |\nabla\ueps|^2
	\qquad \mbox{for all $t>0$ and } \eps\in (0,1),
  \eas
  and since herein
  \bas
	\int_{\{\ueps\ge 1\}} \ueps^{2m-4}|\nabla\ueps|^2
	&=& \int_{\{\ueps\ge 1\}} \ueps^{m-1} \cdot \ueps^{m-3} |\nabla\ueps|^2 \\
	&\le& 2^{3-m} \io (\ueps+\eps)^{m-1} (\ueps+1)^{m-3} |\nabla\ueps|^2
	\qquad \mbox{for all $t>0$ and } \eps\in (0,1)
  \eas
  if $m\le 2$ and, clearly,
  \bas
	\int_{\{\ueps\ge 1\}} \ueps^{2m-4}|\nabla\ueps|^2
	&\le& (\ueps+\eps)^{2m-4} |\nabla\ueps|^2
	\qquad \mbox{for all $t>0$ and } \eps\in (0,1)
  \eas
  if $m>2$, by combining Lemma \ref{lem1} with Lemma \ref{lem43} we obtain functions $k_i: (0,\infty) \to (0,\infty)$, $i\in \{1,2\}$,
  for which (K) holda and which are such that when $T>0$,
  \be{44.4}
	\|P(\ueps)\|_{L^\frac{1}{m-1}(\Omega)} \le k_1(T)
	\qquad \mbox{for all $t\in (0,T+1)$ and } \eps\in (0,1)
  \ee
  and 
  \be{44.5}
	\int_t^{t+1} \|\nabla P(\ueps(\cdot,s))\|_{L^2(\Omega)}^2 ds
	\le k_2(T)
	\qquad \mbox{for all $t\in (0,T)$ and } \eps\in (0,1).
  \ee
  As the Gagliardo-Nirenberg inequality provides $c_1>0$ fulfilling
  \bas
	\io |\varphi|^\frac{2m-1}{m-1} \le c_1 \|\nabla\varphi\|_{L^2(\Omega)}^2 \|\varphi\|_{L^\frac{1}{m-1}(\Omega)}^\frac{1}{m-1}
	+ c_1\|\varphi\|_{L^\frac{1}{m-1}(\Omega)}^\frac{2m-1}{m-1}
	\qquad \mbox{for all } \varphi\in W^{1,2}(\Omega),
  \eas
  we thus infer that for any $T>0$,
  \bas
	\int_t^{t+1} \io P^\frac{2m-1}{m-1}(\ueps)
	&\le& c_1 \int_t^{t+1} \|\nabla P(\ueps(\cdot,s))\|_{L^2(\Omega)}^2 
		\|P(\ueps(\cdot,s))\|_{L^\frac{1}{m-1}(\Omega)}^\frac{1}{m-1} ds \\
	& & + c_1 \int_t^{t+1} \|P(\ueps(\cdot,s))\|_{L^\frac{1}{m-1}(\Omega)}^\frac{2m-1}{m-1} ds \\
	&\le& c_1 k_1^\frac{1}{m-1}(T) k_2(T)
	+ c_1 k_1^\frac{2m-1}{m-1}(T)
	\qquad \mbox{for all $t\in (0,T)$ and } \eps\in (0,1),
  \eas
  from which (\ref{44.1}) immediately follows thanks to (\ref{44.3}) and the trivial fact that $\ueps^{2m-1} \le 3^{2m-1}$
  in $\{\ueps\le 3\}$.
\qed
\mysection{Estimating $\|\veps\|_{W^{1,q}(\Omega)}$ for some $q>2$ when $m>\frac{3}{2}$}
%
%
%
%
%
%
%
%
Now an observation of crucial importance to our approach asserts a bound for $\veps$ with respect to the norm in $W^{1,q}(\Omega)$
with some $q>2$, provided that the integrability exponent in Lemma \ref{lem44} can be chosen to be superquadratic.
This circumstance can actually be viewed as the core of our requirement on $m$ in Theorem \ref{theo68} and Theorem \ref{theo69}.
\begin{lem}\label{lem5}
  Let $m>\frac{3}{2}$. Then one can find a function $K: (0,\infty) \to (0,\infty)$ such that (K) holds, and that if $T>0$ then
  \be{5.1}
	\|\veps(\cdot,t)\|_{W^{1,2m-1}(\Omega)} \le K(T)
	\qquad \mbox{for all $t\in (0,T)$ and } \eps\in (0,1).
  \ee
\end{lem}
\proof
  As $p:=2m-1$ satisfies $p>2$, a Gagliardo-Nirenberg interpolation corresponding to the continuous embeddings
  $W^{1,p}(\Omega) \hra L^\infty(\Omega) \hra L^1(\Omega)$ warrants the existence of $c_1>0$ such that
  \be{5.02}
	\|\varphi\|_{L^\infty(\Omega)} \le c_1 \|\varphi\|_{W^{1,p}(\Omega)}^a \|\varphi\|_{L^1(\Omega)}^{1-a}
	\qquad \mbox{for all } \varphi\in W^{1,p}(\Omega),
  \ee
  with the number $a:=\frac{2p}{3p-2} \in (0,1)$ satisfying
  \bas
	\frac{1}{p} + \frac{p-1}{pa}
	= \frac{1}{p} + \frac{(p-1)(3p-2)}{2p^2}
	= \frac{3p^2 -3p+2}{2p^2}
	>\frac{1}{2},
  \eas
  because $2p^2-3p+2=2(p-1)^2+p>0$.
  We can therefore pick $q\in (1,2)$ sufficiently close to $2$ such that
  \be{5.2}
	\frac{1}{q} < \frac{1}{p} + \frac{p-1}{pa},
  \ee
  and thereupon invoke known smoothing properties of the Neumann heat semigroup $(e^{\sigma\Delta})_{\sigma\ge 0}$ on $\Omega$
  (\cite{win_JDE2010}) to fix positive constants $c_2, c_3$ and $c_4$ fulfilling
  \be{5.3}
	\|e^{\sigma\Delta}\varphi\|_{W^{1,p}(\Omega)} 
	\le c_2 \sigma^{-\alpha} \|\varphi\|_{W^{1,q}(\Omega)}
	\qquad \mbox{for all $\sigma\in (0,2)$ and } \varphi\in C^1(\bom)
  \ee
  and
  \be{5.33}
	\|e^{\sigma\Delta}\varphi\|_{W^{1,p}(\Omega)} \le c_3 \|\varphi\|_{W^{1,\infty}(\Omega)}
	\qquad \mbox{for all $\sigma\in (0,2)$ and } \varphi\in C^1(\bom)
  \ee
  as well as
  \be{5.4}
	\|e^{\sigma\Delta} \varphi\|_{W^{1,p}(\Omega)} \le c_4 \sigma^{-\frac{1}{2}} \|\varphi\|_{L^p(\Omega)}
	\qquad \mbox{for all $\sigma\in (0,2)$ and } \varphi\in C^0(\bom),
  \ee
  where $\alpha:=\frac{1}{q}-\frac{1}{p}>0$.
  Apart from that, Lemma \ref{lem1} together with Lemma \ref{lem44}, (\ref{B}), Lemma \ref{lem3} and (\ref{init}) provides functions
  $k_i: (0,\infty) \to (0,\infty)$, $i\in \{1,2,3,4\}$, which satisfy (K) and are such that for all $T>0$,
  \be{5.5}
	\|\veps(\cdot,t)\|_{L^1(\Omega)} \le k_1(T)
	\qquad \mbox{for all $t\in (0,T)$ and } \eps\in (0,1)
  \ee
  and 
  \be{5.6}
	\int_{(t-2)_+}^t \io \ueps^p \le k_2(T)
	\qquad \mbox{for all $t\in (0,T)$ and } \eps\in (0,1)
  \ee
  as well as
  \be{5.7}
	\|B_2(\cdot,t)\|_{L^p(\Omega)} \le k_3(T)
	\qquad \mbox{for all $t\in (0,T)$ and } \eps\in (0,1),
  \ee
  and that for any such $T$, each $t_0\in (0,T)$ and arbitrary $\eps\in (0,1)$ we can find $\ts\ge 0$
  $t_\star=t_\star(t_0,\eps) \ge 0$ with the properties that
  \be{5.8}
	\ts \in \big((t_0-2)_+, (t_0-1)_+\big)
	\quad \mbox{and} \quad
	\|\veps(\cdot,t_\star)\|_{W^{1,q}(\Omega)} \le k_4(T)
	\qquad \mbox{if } t_0\ge 2,
  \ee
  and that
  \be{5.87}
	\ts=0
	\quad \mbox{and} \quad
	\|\veps(\cdot,t_\star)\|_{W^{1,\infty}(\Omega)} \le c_5:=\|v_0\|_{W^{1,\infty}(\Omega)}
	\qquad \mbox{if } t_0\in (0,2).
  \ee
  Now given $T>0, t_0\in (0,T)$ and $\eps\in (0,1)$, taking $t_\star=t_\star(t_0,\eps)$ as thus specified we estimate the number
  \be{5.89}
	M:=\left\{ \begin{array}{ll}
	\displaystyle
	\sup\limits_{t\in (t_\star,t_0]} \Big\{ (t-t_\star)^\alpha \|\veps(\cdot,t)\|_{W^{1,p}(\Omega)} \Big\}
	\qquad \mbox{if } t_0\ge 2, \\[2mm]
	\displaystyle
	\sup\limits_{t\in (t_\star,t_0]} \|\veps(\cdot,t)\|_{W^{1,p}(\Omega)}
	\end{array} \right.
  \ee
  by relying on a Duhamel representation associated with the second sub-problem of (\ref{0eps}) to see that due to (\ref{5.4}),
  \bea{5.9}
	\|\veps(\cdot,t)\|_{W^{1,p}(\Omega)}
	&=& \bigg\| e^{(t-t_\star)(\Delta-1)} \veps(\cdot,t_\star)
	+ \int_{t_\star}^t e^{(t-s)(\Delta-1)} \frac{\ueps(\cdot,s)\veps(\cdot,s)}{1+\eps\ueps(\cdot,s)\veps(\cdot,s)} ds \nn\\
	& & \hspace*{35mm}
	+ \int_{t_\star}^t e^{(t-s)(\Delta-1)} B_2(\cdot,s) ds \bigg\|_{W^{1,p}(\Omega)} \nn\\
	&\le& e^{-(t-t_\star)} \|e^{(t-\ts)\Delta} \veps(\cdot,\ts)\|_{W^{1,p}(\Omega)} \nn\\
	& & + c_4 \int_{\ts}^t (t-s)^{-\frac{1}{2}} e^{-(t-s)} 
		\bigg\| \frac{\ueps(\cdot,s)\veps(\cdot,s)}{1+\eps\ueps(\cdot,s)\veps(\cdot,s)} \bigg\|_{L^p(\Omega)} ds \nn\\
	& & + c_4 \int_{\ts}^t (t-s)^{-\frac{1}{2}} e^{-(t-s)} \|B_2(\cdot,s)\|_{L^p(\Omega)} ds
	\qquad \mbox{for all $t\in (\ts,t_0]$,}
  \eea
  because $t_0-\ts \in (0,2)$.
  Here if $t_0\ge 2$, then by (\ref{5.3}) and (\ref{5.8}),
  \be{5.10}
	e^{-(t-\ts)} \|e^{(t-\ts)\Delta} \veps(\cdot,\ts)\|_{W^{1,p}(\Omega)}
	\le c_2 (t-\ts)^{-\alpha} \|\veps(\cdot,\ts)\|_{W^{1,q}(\Omega)}
	\le c_2 k_4(T) (t-\ts)^{-\alpha}
	\qquad \mbox{for all } t\in (\ts,t_0],
  \ee
  and if $t_0<2$, then by (\ref{5.33}) and (\ref{5.87}),
  \be{5.100}
	e^{-(t-\ts)} \|e^{(t-\ts)\Delta} \veps(\cdot,\ts)\|_{W^{1,p}(\Omega)}
	\le c_3\|\veps(\cdot,\ts)\|_{W^{1,\infty}(\Omega)} \le c_3 c_5,
  \ee	
  while (\ref{5.7}) asserts that
  \bea{5.11}
	c_4 \int_{\ts}^t (t-s)^{-\frac{1}{2}} e^{-(t-s)} \|B_2(\cdot,s)\|_{L^p(\Omega)} ds
	&\le& c_4 k_3(T) \int_{\ts}^t (t-s)^{-\frac{1}{2}} ds \nn\\
	&=& 2 c_4 k_3(T) (t-\ts)^\frac{1}{2} \nn\\
	&\le& 2^\frac{3}{2} c_4 k_3(T)
	\qquad \mbox{for all } t\in (\ts,t_0].
  \eea
  In order to appropriately cope with the crucial second last summand on the right of (\ref{5.9}), we first concentrate on the case when
  $t_0\ge 2$, in which we apply (\ref{5.02}) together with (\ref{5.5}) and recall our respective definition of $M$ from (\ref{5.89})
  to find that
  \bas
	\bigg\| \frac{\ueps(\cdot,s)\veps(\cdot,s)}{1+\eps\ueps(\cdot,s)\veps(\cdot,s)} \bigg\|_{L^p(\Omega)}
	&\le& \|\ueps(\cdot,s)\|_{L^p(\Omega)} \|\veps(\cdot,s)\|_{L^\infty(\Omega)} \\
	&\le& c_1 \|\ueps(\cdot,s)\|_{L^p(\Omega)} \|\veps(\cdot,s)\|_{W^{1,p}(\Omega)}^a \|\veps(\cdot,s)\|_{L^1(\Omega)}^{1-a} \\
	&\le& c_1 k_1^{1-a}(T) M^a \|\ueps(\cdot,s)\|_{L^p(\Omega)} (s-\ts)^{-\alpha a}
	\qquad \mbox{for all } s\in (\ts,t_0),
  \eas
  so that by the H\"older inequality,
  \bea{5.12}
	& & \hspace*{-15mm}
	c_4 \int_{\ts}^t (t-s)^{-\frac{1}{2}} e^{-(t-s)} 
		\bigg\| \frac{\ueps(\cdot,s)\veps(\cdot,s)}{1+\eps\ueps(\cdot,s)\veps(\cdot,s)} \bigg\|_{L^p(\Omega)} ds \nn\\
	& & \hspace*{-8mm}
	\le \ c_1 c_4 k_1^{1-a}(T) M^a \int_{\ts}^t (t-s)^{-\frac{1}{2}} (s-\ts)^{-\alpha a} \|\ueps(\cdot,s)\|_{L^p(\Omega)} ds \nn\\
	& & \hspace*{-8mm}
	\le \ c_1 c_4 k_1^{1-a}(T) M^a \cdot \bigg\{ \int_{\ts}^t \io \ueps^p \bigg\}^\frac{1}{p} \cdot
		\bigg\{ \int_{\ts}^t (t-s)^{-\frac{p}{2(p-1)}} (s-\ts)^{-\frac{p\alpha a}{p-1}} ds \bigg\}^\frac{p-1}{p}
	\ \mbox{for all } t\in (\ts,t_0].
  \eea
  Here,
  \bas
	\int_{\ts}^t (t-s)^{-\frac{p}{2(p-1)}} (s-\ts)^{-\frac{p\alpha a}{p-1}} ds
	= c_6 (t-\ts)^{1-\frac{p}{2(p-1)} - \frac{p\alpha a}{p-1}}
	\qquad \mbox{for all } t>\ts,
  \eas
  with $c_6:=\int_0^1 (1-\sigma)^{-\frac{p}{2(p-1)}} \sigma^{-\frac{p\alpha a}{p-1}} d\sigma$ being finite,
  because the inequalities $p>2, q>1$ and $a<1$ imply that $\frac{p}{2(p-1)}<1$ and $\frac{p\alpha a}{p-1} = \frac{(p-q)a}{(p-1)q}<a<1$.
  In view of (\ref{5.6}), from (\ref{5.12}) we therefore obtain that
  \bea{5.122}
	& & \hspace*{-20mm}
	c_4 \int_{\ts}^t (t-s)^{-\frac{1}{2}} e^{-(t-s)} 
		\bigg\| \frac{\ueps(\cdot,s)\veps(\cdot,s)}{1+\eps\ueps(\cdot,s)\veps(\cdot,s)} \bigg\|_{L^p(\Omega)} ds \nn\\
	&\le& c_1 c_4 c_6^\frac{p-1}{p} k_1^{1-a}(T) k_2^\frac{1}{p}(T) M^a (t-\ts)^{\frac{p-1}{p}-\frac{1}{2}-\alpha a}
	\qquad \mbox{for all } t\in (\ts,t_0],
  \eea
  which combined with (\ref{5.10}) and (\ref{5.11}) shows that (\ref{5.9}) implies the inequality
  \bea{5.13}
	(t-\ts)^\alpha \|\veps(\cdot,t)\|_{W^{1,p}(\Omega)}
	&\le& c_2 k_4(T)
	+ c_1 c_4 c_6^\frac{p-1}{p} k_1^{1-a}(T) k_2^\frac{1}{p}(T) M^a (t-\ts)^{\frac{p-1}{p} - \frac{1}{2} - \alpha a + \alpha} 
		\nn\\[2mm]
	& & + 2^\frac{3}{2} c_4 k_3(T) (t-\ts)^\alpha \nn\\[2mm]
	&\le& k_5(T) + k_5(T) M^a
	\qquad \mbox{for all } t\in (\ts,t_0]
  \eea
  with
  \bas
	k_5(T):=\max \Big\{ c_2 k_4(T) + 2^\frac{3}{2} c_4 k_3(T) \cdot 2^\alpha \, , \, 
		c_1 c_4 c_6^\frac{p-1}{p} k_1^{1-a}(T) k_2^\frac{1}{p}(T) \cdot 2^{\frac{p-1}{p}-\frac{1}{2}-\alpha a + \alpha} \Big\},
  \eas
  because again since $p>2$ and $a<1$,
  \bas
	\frac{p-1}{p} - \frac{1}{2} - \alpha a + \alpha
	= \frac{p-2}{2p} + (1-a)\alpha >0.
  \eas
  As a further consequence of the fact that $a<1$, (\ref{5.13}) finally entails that
  \bas
	M \le k_6(T):=\max \Big\{ 1 \, , \, (2k_5(T))^\frac{1}{1-a} \Big\},
  \eas
  from which by the definition of $M$ in (\ref{5.89}) it particularly follows that whenever $t_0\ge 2$,
  \be{5.14}
	\|\veps(\cdot,t_0)\|_{W^{1,p}(\Omega)} 
	\le (t_0-\ts)^{-\alpha} M \le k_6(T),
  \ee
  because $t_0-\ts\ge 1$ and $\alpha\ge 0$.\abs
  If $t_0\in (0,2)$, however, then referring to the respective part in (\ref{5.89}) enables us to actually simplify the above reasoning
  so as to infer, in a way similar to that in (\ref{5.12}) and (\ref{5.122}), that
  \bas
	& & \hspace*{-20mm}
	c_4 \int_{\ts}^t (t-s)^{-\frac{1}{2}} e^{-(t-s)} 
		\bigg\| \frac{\ueps(\cdot,s)\veps(\cdot,s)}{1+\eps\ueps(\cdot,s)\veps(\cdot,s)} \bigg\|_{L^p(\Omega)} ds \nn\\
	&\le& c_1 c_4 k_1^{1-a}(T) M^a \cdot \bigg\{ \int_0^t \io \ueps^p \bigg\}^\frac{1}{p} \cdot
		\bigg\{ \int_0^t (t-s)^{-\frac{p}{2(p-1)}} ds \bigg\}^\frac{p-1}{p} \\
	&\le& c_1 c_4 c_7^\frac{p-1}{p} k_1^{1-a}(T) k_2^\frac{1}{p}(T) M^a
	\qquad \mbox{for all } t\in (\ts,t_0] \equiv (0,t_0],
  \eas
  with $c_7:=\int_0^2 \sigma^{-\frac{p}{2(p-1)}} d\sigma \equiv \frac{2(p-1)}{p-2} \cdot 2^\frac{p-2}{2(p-1)}$.
  In this case now relying on (\ref{5.100}) instead of (\ref{5.10}), from (\ref{5.9}) and (\ref{5.11}) we thus infer that
  \bas
	\|\veps(\cdot,t)\|_{W^{1,p}(\Omega)}
	\le c_3 c_5 + c_1 c_4 c_7^\frac{p-1}{p} k_1^{1-a}(T) k_2^\frac{1}{p}(T) M^a
	+ 2^\frac{3}{2} c_4 k_3(T)
	\qquad \mbox{for all } t\in (0,t_0]
  \eas
  and that hence
  \bas
	M \le k_7(T) + k_7(T) M^a,
  \eas
  where $k_7(T):=\max\{ c_3 c_5 + 2^\frac{3}{2} c_4 k_3(T) \, , \, c_1 c_4 c_7^\frac{p-1}{p} k_1^{1-a}(T) k_2^\frac{1}{p}(T) \}$.
  Again since $a<1$, this especially shows that for any such $t_0$,
  \bas
	\|\veps(\cdot,t_0)\|_{W^{1,p}(\Omega)} \le M \le \max \Big\{ 1 \, , \, (2k_7(T))^\frac{1}{1-a} \Big\},
  \eas
  which together with (\ref{5.14}) yields the claimed conclusion.
\qed
\mysection{Boundedness properties in $L^\infty(\Omega)\times W^{1,q}(\Omega)$ for arbitrary $q>2$}
With the knowledge from Lemma \ref{lem5} at hand, we can successively improve our information about regularity in the course
of a three-step bootstrap procedure, the first part of which is concerned with bounds on $\ueps$ in $L^p(\Omega)$
for arbitrarily large finite $p$.
\begin{lem}\label{lem55}
  Let $m>\frac{3}{2}$ and $p>\max\{1 \, , m-1+\frac{2m-3}{2m-1}\}$. 
  Then there exists $K \equiv K^{(p)}: (0,\infty) \to (0,\infty)$ such that (K) is valid and that whenever $T>0$,
  \be{55.1}
	\io \ueps^p(\cdot,t) \le K(T)
	\qquad \mbox{for all $t\in (0,T)$ and } \eps\in (0,1).
  \ee
\end{lem}
\proof
  On testing the first equation in (\ref{0eps}) by $\ueps^{p-1}$ and using Young's inequality, we see that
  \bea{55.2}
	& & \hspace*{-20mm}
	\frac{1}{p} \frac{d}{dt} \io \ueps^p
	+ \frac{2(p-1)}{(m+p-1)^2} \io |\nabla\ueps^\frac{m+p-1}{2}|^2 \nn\\
	&=& \frac{p-1}{2} \io \ueps^{m+p-3} |\nabla \ueps|^2
	- (p-1) \io \ueps^{p-2} (\ueps+\eps)^{m-1} |\nabla\ueps|^2 \nn\\
	& & + (p-1)\chi \io \ueps^{p-1} \nabla\ueps \cdot \frac{\nabla\veps}{\veps}
	- \io \ueps^p \veps 
	+ \io \ueps^{p-1} B_1 \nn\\
	&\le& \frac{(p-1)\chi^2}{2} \io \ueps^{-m+p-1} \frac{|\nabla\veps|^2}{\veps^2}
	- \io \ueps^p \veps
	+ \io \ueps^{p-1} B_1
	\qquad \mbox{for all } t>0,
  \eea
  where using Lemma \ref{lem01} along with (\ref{B}) and again Young's inequality we can find $k_i: (0,\infty) \to (0,\infty)$,
  $i\in \{1,2,3\}$, fulfilling (K) and such that for all $T>0$,
  \bea{55.3}
	- \io \ueps^p \veps
	+ \io \ueps^{p-1} B_1
	&\le& - k_1(T) \io \ueps^p
	+ k_2(T) \io \ueps^{p-1} \nn\\
	&\le& - \frac{k_1(T)}{2} \io \ueps^p
	+ k_3(T)
	\qquad \mbox{for all $t\in (0,T)$ and } \eps\in (0,1).
  \eea
  Apart from that, Lemma \ref{lem5} in conjunction with Lemma \ref{lem01} entails the existence of $k_4: (0,\infty) \to (0,\infty)$
  such that (K) holds, and that if $T>0$ then
  \bas
	\io \frac{|\nabla\veps|^{2m-1}}{\veps^{2m-1}} \le k_4(T)
	\qquad \mbox{for all $t\in (0,T)$ and } \eps\in (0,1),
  \eas
  whence utilizing the H\"older inequality we find that for any such $T$,
  \bea{55.4}
	& & \hspace*{-30mm}
	\frac{(p-1)\chi^2}{2} \io \ueps^{-m+p-1} \frac{|\nabla\veps|^2}{\veps^2} \nn\\
	&\le& \frac{(p-1)\chi^2}{2} \cdot \bigg\{ \io \ueps^\frac{(2m-1)(-m+p-1)}{2m-3} \bigg\}^\frac{2m-3}{2m-1} \cdot
		\bigg\{ \io \frac{|\nabla \veps|^{2m-1}}{\veps^{2m-1}} \bigg\}^\frac{2m-3}{2m-1} \nn\\
	&\le& k_5(T) \cdot \bigg\{ \io \ueps^\frac{(2m-1)(-m+p-1)}{2m-3} \bigg\}^\frac{2m-3}{2m-1}
	\qquad \mbox{for all $t\in (0,T)$ and } \eps\in (0,1)
  \eea
  with $k_5(T):=\frac{(p-1)\chi^2}{2} k_4^\frac{2m-3}{2m-1}(T)$.
  Now since $\frac{2}{m+p-1}<\frac{2}{m+p-1} \cdot \frac{(2m-1)(-m+p+1)}{2m-3}$ due to the fact that
  $-m+p+1>\frac{2m-3}{2m-1}$ by assumption on $p$, the Gagliardo-Nirenberg inequality applies so as to say that with
  \be{55.44}
	a:=\frac{(2m-1)(-m+p+1)-2m+3}{(2m-1)(-m+p+1)} \, \in (0,1)
  \ee
  and some $c_1=c_1(p)>0$ we have
  \bea{55.5}
	\bigg\{ \io \ueps^\frac{(2m-1)(-m+p-1)}{2m-3} \bigg\}^\frac{2m-3}{2m-1}
	&=& \|\ueps^\frac{m+p-1}{2}\|_{L^{\frac{2}{m+p-1} \cdot \frac{2m-1)(-m+p+1)}{2m-3}}(\Omega)}^\frac{2(-m+p+1)}{m+p-1} \nn\\
	&\le& c_1\|\nabla \ueps^\frac{m+p-1}{2}\|_{L^2(\Omega)}^\frac{2(-m+p+1)a}{m+p-1}
		\|\ueps^\frac{m+p-1}{2}\|_{L^\frac{2}{m+p-1}(\Omega)}^\frac{2(-m+p+1)(1-a)}{m+p-1} \nn\\
	& & + \|\ueps^\frac{m+p-1}{2}\|_{L^\frac{2}{m+p-1}(\Omega)}^\frac{2(-m+p+1)}{m+p-1}
	\qquad \mbox{for all $t>0$ and } \eps\in (0,1).
  \eea
  Here we recall that Lemma \ref{lem1} provides $k_6: (0,\infty)\to (0,\infty)$ such that (K) holds and that for all $T>0$,
  \bas
	\|\ueps^\frac{m+p-1}{2}\|_{L^\frac{2}{m+p-1}(\Omega)}^\frac{2}{m+p-1}
	= \io \ueps \le k_6T)
	\qquad \mbox{for all $t\in (0,T)$ and } \eps\in (0,1),
  \eas
  and furthermore we note that according to (\ref{55.44}),
  \bas
	\frac{(-m+p+1)a}{m+p-1} -1
	&=& \frac{(2m-1)(-m+p+1) -2m+3}{(2m-1)(m+p-1)} -1 \\
	&=& \frac{(2m-1)(-2m+2) -2m+3}{(2m-1)(m+p-1)} \\
	&=& -\frac{2(m-1)}{m+p-1} - \frac{2m-3}{(2m-1)(m+p-1)} \\[2mm]
	&<& 0,
  \eas
  so that $\theta:=\frac{m+p-1}{(-m+p+1)a}$ satisfies $\theta >1$.
  An application of Young's inequality to (\ref{55.5}) therefore yields functions $k_i: (0,\infty) \to (0,\infty)$, $i\in \{7,8\}$,
  for which (K) is valid, and which are such that for all $T>0$,
  \bas
	& & \hspace*{-30mm}
	k_5(T) \cdot \bigg\{ \io \ueps^\frac{(2m-1)(-m+p-1)}{2m-3} \bigg\}^\frac{2m-3}{2m-1} \nn\\
	&\le& k_7(T) \|\nabla\ueps^\frac{m+p-1}{2}\|_{L^2(\Omega)}^\frac{2(-m+p+1)a}{m+p-1} + k_7(T) \\
	&\le& \frac{2(p-1)}{(m+p-1)^2} \io |\nabla\ueps^\frac{m+p-1}{2}|^2
	+ k_8(T)
	\qquad \mbox{for all $t\in (0,T)$ and } \eps\in (0,1).
  \eas 
  Together with (\ref{55.3}) and (\ref{55.4}) inserted into (\ref{55.2}), this shows that for each $T>0$ we have
  \bas
	\frac{1}{p} \frac{d}{dt} \io \ueps^p 
	+ \frac{k_1(T)}{2} \io \ueps^p
	\le k_3(T) + k_8(T)
	\qquad \mbox{for all $t\in (0,T)$ and } \eps\in (0,1),
  \eas
  which results in (\ref{55.1}) by means of an evident ODE comparison argument.
\qed
This in turn improves our knowledge on the second solution component:
\begin{lem}\label{lem56}
  Let $m>\frac{3}{2}$ and $q>2$.
  Then one can find $K \equiv K^{(q)}: (0,\infty) \to (0,\infty)$ such that (K) holds, and that given any $T>0$ we have
  \be{56.1}
	\|\veps(\cdot,t)\|_{W^{1,q}(\Omega)} \le K(T)
	\qquad \mbox{for all $t\in (0,T)$ and } \eps\in (0,1).
  \ee
\end{lem}
\proof
  As $W^{1,2m-1}(\Omega) \hra L^\infty(\Omega)$ due to the hypothesis $m>\frac{3}{2}$, Lemma \ref{lem5} together with Lemma \ref{lem55}
  and (\ref{B}) in particular yields $k_1: (0,\infty) \to (0,\infty)$ fulfilling (K) and such 
  that writing $f_\eps(x,t):=\frac{\ueps}{1+\eps\ueps\veps}(x,t)+B_2(x,t)$, $(x,t)\in \Omega\times (0,\infty)$, $\eps\in (0,1)$,
  for all $T>0$ we have
  \bas
	\|f_\eps(\cdot,t)\|_{L^p(\Omega)} \le k_1(T)
	\qquad \mbox{for all $t\in (0,T)$ and } \eps\in (0,1).
  \eas
  Therefore, (\ref{56.1}) can be derived by straightforward application of well-known regularization estimates for the
  Neumann heat semigroup (\cite{win_JDE2010}) to the inhomogeneous linear heat equation $v_{\eps t}=\Delta \veps + f_\eps$.
\qed
When combined with Lemma \ref{lem55}, through a standard argument 
the latter in fact asserts a boundedness feature of $\ueps$ even with respect to the norm in $L^\infty(\Omega)$.
\begin{lem}\label{lem6}
  Let $m>\frac{3}{2}$.
  Then there exists $K: (0,\infty) \to (0,\infty)$ fulfilling (K) such that for all $T>0$,
  \bas
	\|\ueps(\cdot,t)\|_{L^\infty(\Omega)} \le K(T)
	\qquad \mbox{for all $t\in (0,T)$ and } \eps\in (0,1).
  \eas
\end{lem}
\proof
  This can readily be obtained from the bounds provided by Lemma \ref{lem55} and Lemma \ref{lem56} through a standard application
  of a Moser-type recursive argument (cf.~e.g.~\cite[Lemma A.1]{taowin_subcritical}).
\qed
\mysection{Further compactness properties and regularity in time}
For our mere existence statement in Theorem \ref{theo68}, tracking a possible dependence of estimates on the asymptotic behavior
of $B_1$ and $B_2$ seems unnecessary;
the next three statements preparing our limit procedure $\eps\searrow 0$ will therefore not involve our hypothesis (K),
but rather exclusively provide information on arbitrary but fixed time intervals. 
Our first observation in this regard is an essentially immediate consequence Lemma \ref{lem43} when combined with the
boundedness information from Lemma \ref{lem6}.
\begin{lem}\label{lem51}
  Let $m>\frac{3}{2}$ and
  \be{51.1}
	\alpha \ge \left\{ \begin{array}{ll}
	\displaystyle
	\frac{m+1}{2}
	\qquad & \mbox{if } m\in \Big(\frac{3}{2},2\Big], \\[2mm]
	m-1
	& \mbox{if } m>2.
	\end{array} \right.
  \ee
  Then for all $T>0$ there exists $C(\alpha,T)>0$ such that
  \be{51.2}
	\int_0^T \io |\nabla (\ueps+\eps)^\alpha|^2 \le C(\alpha,T)
	\qquad \mbox{for all } \eps\in (0,1).
  \ee
\end{lem}
\proof
  In view of Lemma \ref{lem6}, given $T>0$ we can fix $c_1(T)>0$ fulfilling
  \be{51.3}
	\ueps \le c_1(T)
	\quad \mbox{in } \Omega\times (0,T)
	\qquad \mbox{for all } \eps\in (0,1).
  \ee
  Therefore, in the case $m\in (\frac{3}{2},2]$ we can use that then (\ref{51.1}) requires that $2\alpha\ge m+1$ to estimate
  \bas
	\frac{1}{\alpha^2} |\nabla (\ueps+\eps)^\alpha|^2
	&=& (\ueps+\eps)^{2\alpha-2} |\nabla\ueps|^2 \\
	&=& \Big\{ (\ueps+\eps)^{m-1} (\ueps+1)^{m-3} |\nabla\ueps|^2 \Big\} \cdot
	(\ueps+\eps)^{2\alpha-m-1} (\ueps+1)^{3-m} \\
	&\le& \Big\{ (\ueps+\eps)^{m-1} (\ueps+1)^{m-3} |\nabla\ueps|^2 \Big\} \cdot
	(c_1(T)+1)^{2\alpha-m-1} (c_1(T)+1)^{3-m}
  \eas
  in $\Omega\times (0,T)$ for all $\eps\in (0,1)$,
  so that in light of Lemma \ref{lem43}, (\ref{51.2}) results upon an integration over $\Omega\times (0,T)$.\abs
  Similarly, if $m>2$ then $2\alpha\ge 2m-2$ by (\ref{51.1}), and thus 
  \bas
	\frac{1}{\alpha^2} |\nabla (\ueps+\eps)^\alpha|^2
	&=& \Big\{ (\ueps+\eps)^{2m-4} |\nabla\ueps|^2 \Big\} \cdot
	(\ueps+\eps)^{2\alpha-2m-2}\\
	&\le& \Big\{ (\ueps+\eps)^{2m-4} |\nabla\ueps|^2 \Big\} \cdot
	(c_1(T)+1)^{2\alpha-2m-2}
	\quad \mbox{in } \Omega\times (0,T)
	\qquad \mbox{for all } \eps\in (0,1),
  \eas
  again implying (\ref{51.2}) due to Lemma \ref{lem43}.
\qed
Now for suitably large $\alpha$, the expressions appearing in (\ref{51.2}) enjoy some favorable time regularity feature:
\begin{lem}\label{lem52}
  Let $m>\frac{3}{2}$ and
  \be{52.1}
	\alpha \ge \left\{ \begin{array}{ll}
	\displaystyle
	2
	\qquad & \mbox{if } m\in \Big(\frac{3}{2},2\Big], \\[2mm]
	\displaystyle
	m-1
	& \mbox{if } m>2.
	\end{array} \right.
  \ee
  Then for all $T>0$ there exists $C(\alpha,T)>0$ such that
  \be{52.2}
	\int_0^T \Big\| \partial_t \Big(\ueps(\cdot,t)+\eps\Big) \Big\|_{(W^{2,2}(\Omega))^\star} dt \le C(\alpha,T)
	\qquad \mbox{for all } \eps\in (0,1).
  \ee
\end{lem}
\proof
  Using (\ref{0eps}), for fixed $t>0$ and $\varphi\in C^\infty(\bom)$ we compute
  \bea{52.3}
	& & \hspace*{-20mm}
	\frac{1}{\alpha} \io \partial_t (\ueps+\eps)^\alpha \varphi \nn\\
	&=& \io (\ueps+\eps)^{\alpha-1} \varphi \cdot \Big\{ \nabla \cdot \big( (\ueps+\eps)^{m-1} \nabla\ueps \big)
	- \chi \nabla \cdot \Big( \frac{\ueps}{\veps} \nabla\veps \Big) - \ueps\veps + B_1 \Big\} \nn\\
	&=& - \io \Big\{ (\alpha-1) (\ueps+\eps)^{\alpha-2} \varphi\nabla\ueps + (\ueps+\eps)^{\alpha-1} \nabla \varphi \Big\} \cdot
		\Big\{ (\ueps+\eps)^{m-1} \nabla\ueps - \chi \frac{\ueps}{\veps}\nabla\veps \Big\} \nn\\
	& & - \io \ueps (\ueps+\eps)^{\alpha-1} \veps\varphi
	+ \io (\ueps+\eps)^{\alpha-1} \varphi \nn\\
	&=& - (\alpha-1) \io (\ueps+\eps)^{m+\alpha-3} |\nabla\ueps|^2 \varphi 
	+ (\alpha-1) \chi \io \ueps (\ueps+\eps)^{\alpha-2} \Big(\nabla\ueps \cdot\frac{\nabla\veps}{\veps}\Big) \varphi \nn\\
	& & - \io (\ueps+\eps)^{m+\alpha-2} \nabla\ueps\cdot\nabla\varphi
	+ \chi \io \ueps(\ueps+\eps)^{\alpha-1} \frac{\nabla\veps}{\veps}\cdot\nabla\varphi \nn\\
	& & - \io \ueps (\ueps+\eps)^{\alpha-1} \veps\varphi
	+ \io (\ueps+\eps)^{\alpha-1} \varphi
	\qquad \mbox{for all } \eps\in (0,1).
  \eea
  Here given $T>0$, we note that Lemma \ref{lem6}, Lemma \ref{lem01}, Lemma \ref{lem5} and (\ref{B}) ensure the existence of positive
  constants $c_i(T)$, $i\in\{1,2,3,4\}$, such that for all $\eps\in (0,1)$,
  \be{52.4}
	\ueps\le c_1(T),
	\quad 
	c_2(T) \le \veps \le c_3(T)
	\quad \mbox{and} \quad
	B_1 \le c_4(T)
	\qquad \mbox{in } \Omega\times (0,T).
  \ee
  Since (\ref{52.1}) especially requires that $\alpha\ge 1$, by using Young's inequality we thus obtain that whenever $t\in (0,T)$ 
  and $\eps\in (0,1)$,
  \bea{52.5}
	\bigg| \chi \io \ueps(\ueps+\eps)^{\alpha-1} \frac{\nabla\veps}{\veps}\cdot\nabla\varphi \bigg|
	&\le& \io \veps^{-\frac{3}{2}} |\nabla \veps|^2
	+ \frac{\chi^2}{4} \io \frac{\ueps^2 (\ueps+\eps)^{2\alpha-2}}{\veps^\frac{1}{2}} |\nabla\varphi|^2 \nn\\
	&\le& \io \veps^{-\frac{3}{2}} |\nabla \veps|^2
	+ \frac{\chi^2}{4} \cdot \frac{c_1^2(T) \cdot (c_1(T)+1)^{2\alpha-2}}{c_2^\frac{1}{2}(T)} \cdot 
		\|\nabla\varphi\|_{L^2(\Omega)}^2
  \eea
  and
  \be{52.6}
	\bigg| - \io \ueps(\ueps+\eps)^{\alpha-1} \veps\varphi \bigg|
	\le c_1(T) \cdot (c_1(T)+1)^{\alpha-1} c_3(T) |\Omega| \cdot \|\varphi\|_{L^\infty(\Omega)}
  \ee
  as well as
  \be{52.7}
	\bigg| \io (\ueps+\eps)^{\alpha-1} B_1\varphi \bigg|
	\le (c_1(T)+1)^{\alpha-1} c_4(T) |\Omega| \cdot \|\varphi\|_{L^\infty(\Omega)}.
  \ee
  Now in the case $m\in (\frac{3}{2},2]$ in which (\ref{52.1}) asserts that 
  $\alpha\ge 2\ge \max\{ \frac{m+1}{2} \, , \, \frac{3-m}{2}\}$,
  the first three summand on the right of (\ref{52.3}) can similarly be estimated according to
  \bea{52.8}
	& & \hspace*{-20mm}
	\bigg| -(\alpha-1) \io (\ueps+\eps)^{m+\alpha-3} |\nabla\ueps|^2 \varphi \bigg| \nn\\
	&\le& (\alpha-1) (c_1(T)+\eps)^{\alpha-2} (c_1(T)+1)^{3-m} \cdot 
		\bigg\{ \io (\ueps+\eps)^{m-1} (\ueps+1)^{m-3} |\nabla\ueps|^2 \bigg\} \cdot \|\varphi\|_{L^\infty(\Omega)}
  \eea
  and
  \bea{52.9}
	& & \hspace*{-20mm}
	\bigg| (\alpha-1)\chi \io \ueps(\ueps+\eps)^{\alpha-2} \Big(\nabla\ueps\cdot\frac{\nabla\veps}{\veps}\Big) \varphi \bigg| \nn\\
	&\le& \veps^{-\frac{3}{2}}|\nabla\veps|^2
	+ \frac{(\alpha-1)^2 \chi^2}{4} \cdot \bigg\{ \io \frac{\ueps^2 (\ueps+\eps)^{2\alpha-4}}{\veps^\frac{1}{2}} |\nabla\ueps|^2
		\bigg\} \cdot \|\varphi\|_{L^\infty(\Omega)} \nn\\
	&\le& \veps^{-\frac{3}{2}}|\nabla\veps|^2 \nn\\
	& & + \frac{(\alpha-1)^2 \chi^2}{4} \cdot \frac{(c_1(T)+\eps)^{-m+2\alpha-1} (c_1(T)+1)^{3-m}}{c_2^\frac{1}{2}(T)} \times \nn\\
	& & \hspace*{50mm} 
	\times \bigg\{ \io (\ueps+\eps)^{m-1} (\ueps+1)^{3-m} |\nabla\ueps|^2 \bigg\} \cdot \|\varphi\|_{L^\infty(\Omega)}
  \eea
  and
  \bea{52.10}
	& & \hspace*{-20mm}
	\bigg| - \io (\ueps+\eps)^{m+\alpha-2} \nabla\ueps\cdot\nabla\varphi \bigg| \nn\\
	&\le& \io (\ueps+\eps)^{2m+2\alpha-4} |\nabla\ueps|^2
	+ \frac{1}{4} \|\nabla\varphi\|_{L^2(\Omega)}^2 \nn\\
	&\le& (c_1(T)+\eps)^{m+2\alpha-3} (c_1(T)+1)^{3-m} \io (\ueps+\eps)^{m-1} (\ueps+1)^{m-3} |\nabla \ueps|^2
	+ \frac{1}{4} \|\nabla\varphi\|_{L^2(\Omega)}^2
  \eea
  for all $t\in (0,T)$ and $\eps\in (0,1)$.
  Since $W^{2,2}(\Omega) \hra L^\infty(\Omega)$, from (\ref{52.3}) and (\ref{52.5})-(\ref{52.10}) we thus infer that 
  if $m\in (\frac{3}{2},2]$ and $\alpha$ satisfies (\ref{52.1}), then for each $T>0$ there exists $c_5(T)>0$ such that
  for all $t\in (0,T)$ and $\eps\in (0,1)$,
  \bas
	\| \partial_t (\ueps+\eps)\|_{(W^{2,2}(\Omega))^\star}
	\le c_5(T) \cdot \bigg\{ \io (\ueps+\eps)^{m-1} (\ueps+1)^{m-3} |\nabla\ueps|^2 
	+ \io \veps^{-\frac{3}{2}} |\nabla\veps|^2 + 1 \bigg\},
  \eas
  so that (\ref{52.2}) results from Lemma \ref{lem43} and Lemma \ref{lem2} upon an integration in time for any such $m$ and $\alpha$.\abs
  If $m>2$, then in view of the accordingly modified form of the estimate in Lemma \ref{lem43}, given $T>0$ we rely
  on the hypothesis $\alpha\ge m-1$ in replacing (\ref{52.8})-(\ref{52.10}) with the inequalities
  \bas
	& & \hspace*{-30mm}
	\bigg| -(\alpha-1) \io (\ueps+\eps)^{m+\alpha-3} |\nabla\ueps|^2 \varphi \bigg| \\
	&\le& (\alpha-1) \cdot (c_1(T)+\eps)^{-m+\alpha+1} \cdot 
		\bigg\{ \io (\ueps+\eps)^{2m-4} |\nabla\ueps|^2 \bigg\} \cdot \|\varphi\|_{L^\infty(\Omega)}
  \eas
  and
  \bas
	& & \hspace*{-20mm}
	\bigg| (\alpha-1)\chi \io \ueps(\ueps+\eps)^{\alpha-2} \Big(\nabla\ueps\cdot\frac{\nabla\veps}{\veps}\Big) \varphi \bigg| \\
	&\le& \veps^{-\frac{3}{2}}|\nabla\veps|^2
	+ \frac{(\alpha-1)^2 \chi^2}{4} \cdot \frac{(c_1(T)+\eps)^{-2m+2\alpha+2}}{c_2^\frac{1}{2}(T)} \cdot
		\bigg\{ \io (\ueps+\eps)^{2m-4} |\nabla\ueps|^2 \bigg\} \cdot \|\varphi\|_{L^\infty(\Omega)}
  \eas
  as well as
  \bas
	\bigg| - \io (\ueps+\eps)^{m+\alpha-2} \nabla\ueps\cdot\nabla\varphi \bigg| 
	\le (c_1(T)+\eps)^{2\alpha} \io (\ueps+\eps)^{2m-4} |\nabla \ueps|^2
	+ \frac{1}{4} \|\nabla\varphi\|_{L^2(\Omega)}^2
  \eas
  for all $t\in (0,T)$ and $\eps\in (0,1)$, and conclude as before.
\qed
Independently from the latter two lemmata, the estimates from Lemma \ref{lem6} and Lemma \ref{lem56}
entail a H\"older regularity property of the second solution component as follows.
\begin{lem}\label{lem666}
  Let $m>\frac{3}{2}$.
  Then for all $T>0$ there exist $\vartheta=\vartheta(T)\in (0,1)$ and $C(T)>0$ such that
  \bas
	\|\veps\|_{C^{\vartheta,\frac{\vartheta}{2}}(\bom\times [0,T])} \le C(T)
	\qquad \mbox{for all } \eps\in (0,1).
  \eas
\end{lem}
\proof
  Once more letting $f_\eps:=\frac{\ueps\veps}{1+\eps\ueps\veps} + B_2$ in $\Omega\times (0,\infty)$ for $\eps\in (0,1)$,
  from Lemma \ref{lem6} and, e.g., Lemma \ref{lem56} we especially know that $(f_\eps)_{\eps\in (0,1)}$ is bounded in
  $L^\infty_{loc}(\bom \times [0,\infty))$. As $v_0$ is H\"older continuous in $\bom$ thanks to (\ref{init}), the claimed estimate
  therefore directly follows from standard theory on H\"older regularity in scalar parabolic equations (\cite{porzio_vespri}).
\qed
\mysection{Passing to the limit. Proof of Theorem \ref{theo68} and Theorem \ref{theo69}}
We are now prepared to construct a solution of (\ref{0}) by means of appropriate compactness arguments, 
where following quite standard precedents, our concept of solvability will be as specified in the following.
\begin{defi}\label{dw}
  Assume that $m\ge 1$, that $\chi\in\R$, and that (\ref{B}) and (\ref{init}) hold. Then a pair $(u,v)$ of functions
  \be{w1}
	\left\{ \begin{array}{l}
	u \in L^m_{loc}(\bom\times [0,\infty))
	\qquad \mbox{and} \\[1mm]
	v \in L^1_{loc}([0,\infty);W^{1,1}(\Omega))
	\end{array} \right.
  \ee
  will be called a {\em global weak solution of (\ref{0})} if $u\ge 0$ and $v>0$ a.e.~in $\Omega\times (0,\infty)$, if
  \be{w2}
	\frac{u}{v} \nabla v
	\mbox{ belongs to } L^1_{loc}(\bom\times [0,\infty);\R^2)
  \ee
  and
  \be{w3}
	uv
	\mbox{ lies in } L^1_{loc}(\bom\times [0,\infty)),
  \ee
  and if for each $\varphi\in C_0^\infty(\bom\times [0,\infty))$ fulfilling $\frac{\partial\varphi}{\partial\nu}=0$
  on $\pO\times (0,\infty)$, and for any $\phi\in C_0^\infty(\bom\times [0,\infty))$, the identities
  \be{w4}
	- \int_0^\infty \io u\varphi_t 
	- \io u_0 \varphi(\cdot,0)
	= \frac{1}{m} \int_0^\infty \io u^m \Delta\varphi
	+ \chi \int_0^\infty \io \frac{u}{v}\nabla v\cdot\nabla\varphi
	- \int_0^\infty \io uv\varphi
	+ \int_0^\infty \io B_1 \varphi
  \ee
  and
  \be{w5}
	- \int_0^\infty \io v\phi_t
	- \io v_0\phi(\cdot,0)
	= - \int_0^\infty \io \nabla v\cdot\nabla \phi
	- \int_0^\infty \io  v\phi
	+ \int_0^\infty \io uv\phi
	+ \int_0^\infty \io B_2 \phi
  \ee
  are valid.
\end{defi}

\begin{figure}[H]
  \center
      \subfloat[$t=0$]{\label{fig:init}\includegraphics[width=0.5\textwidth]{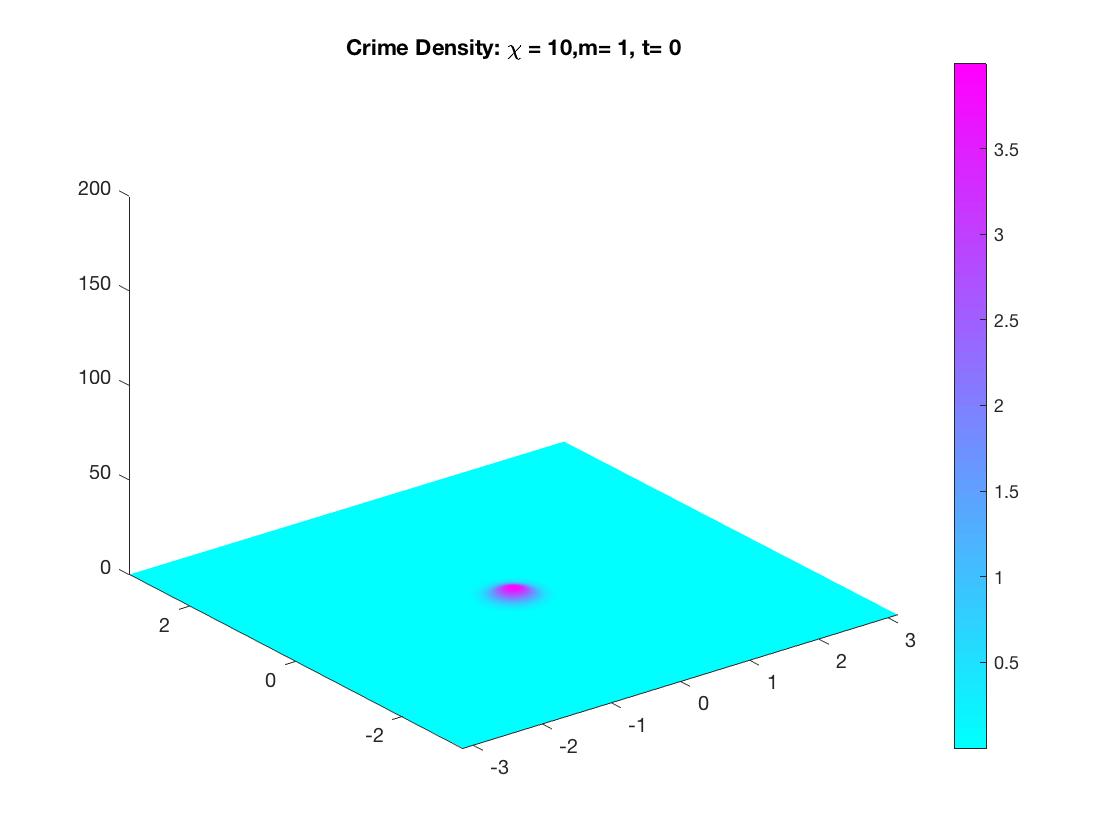}}\;
          \subfloat[$t=.95$]{\label{fig:blowup}\includegraphics[width=0.45\textwidth]{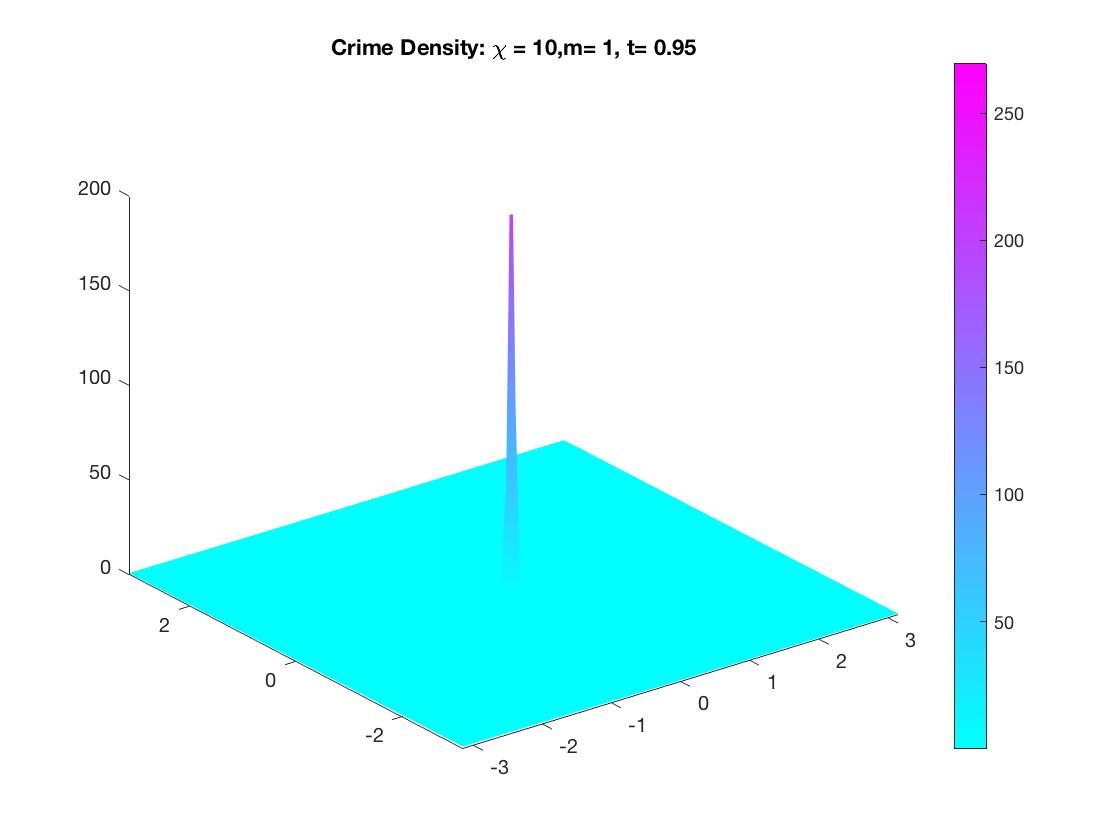}}\;
   \caption{Numerical solutions with $m=1,\;\chi=10$ and $u_0(x) =v_0(x) = \frac{1}{\sqrt{2\pi \sigma^2}}e^{-\frac{|x|^2}{2\sigma^2}}$ with $\sigma = 1/4$. }   \label{fig:short}
\end{figure}

\begin{figure}[H]
  \center
          \subfloat[$t=.95$]{\label{fig:1}\includegraphics[width=0.45\textwidth]{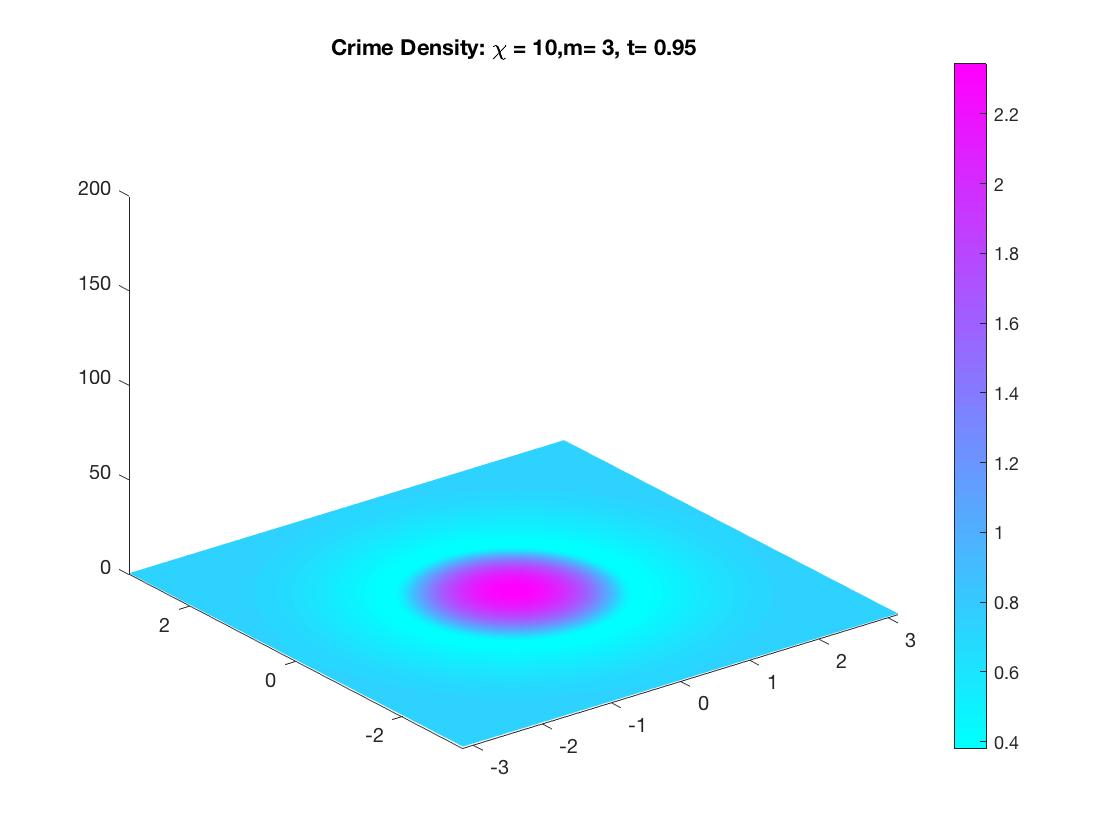}}\;
     \subfloat[$t = 1.2$]{\label{fig:2}\includegraphics[width=0.45\textwidth]{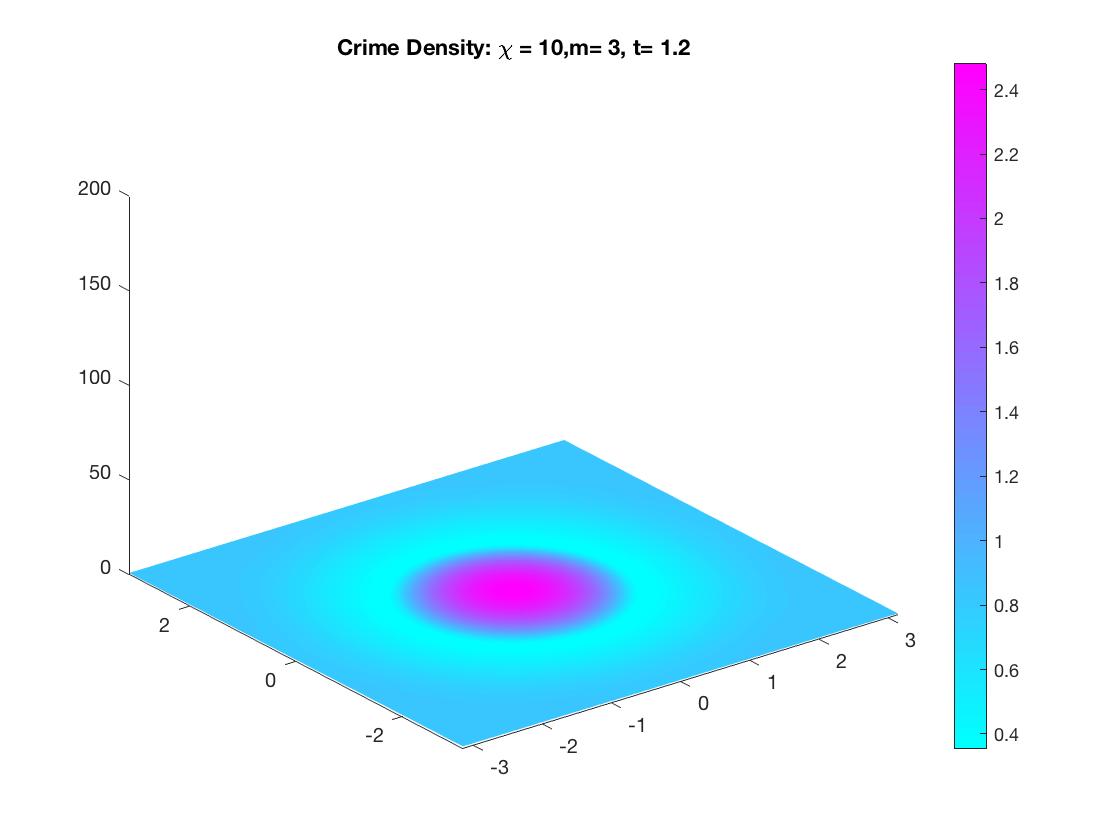}}\\
  \subfloat[$t=1.95$]{\label{fig:3}\includegraphics[width=0.45\textwidth]{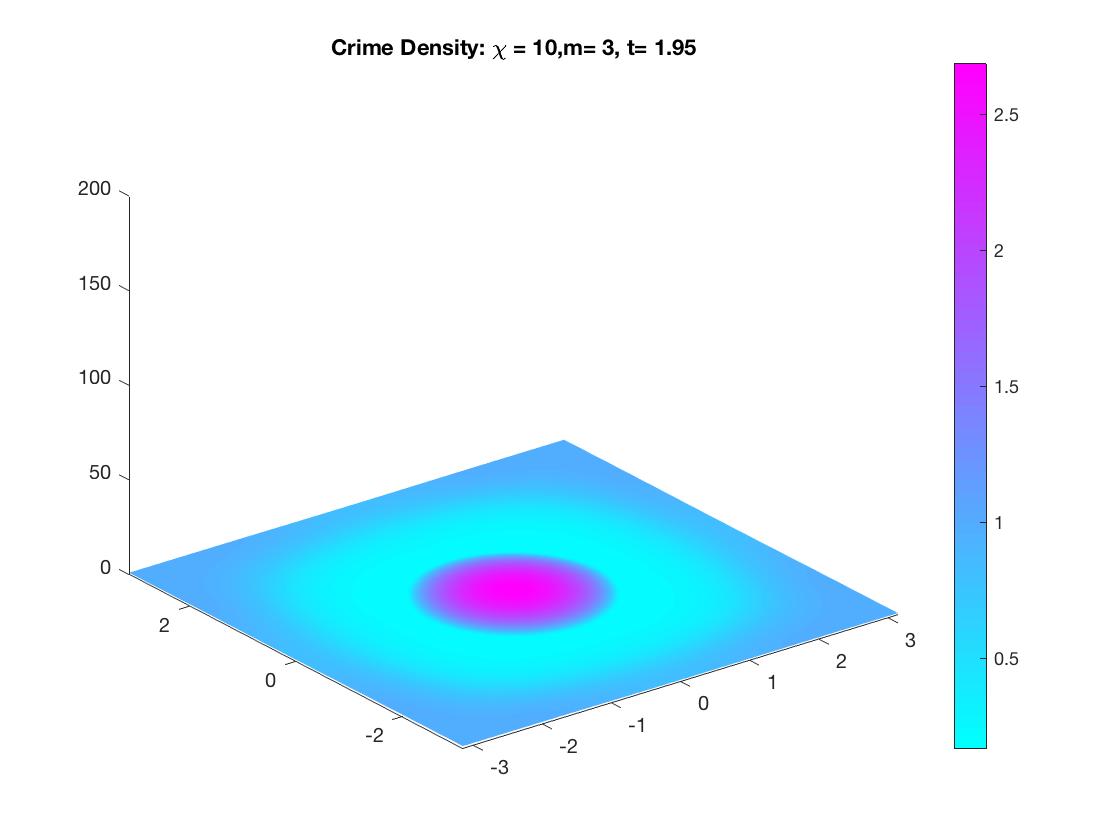}}
  \subfloat[$t = 9.95$]{\label{fig:4}\includegraphics[width=0.45\textwidth]{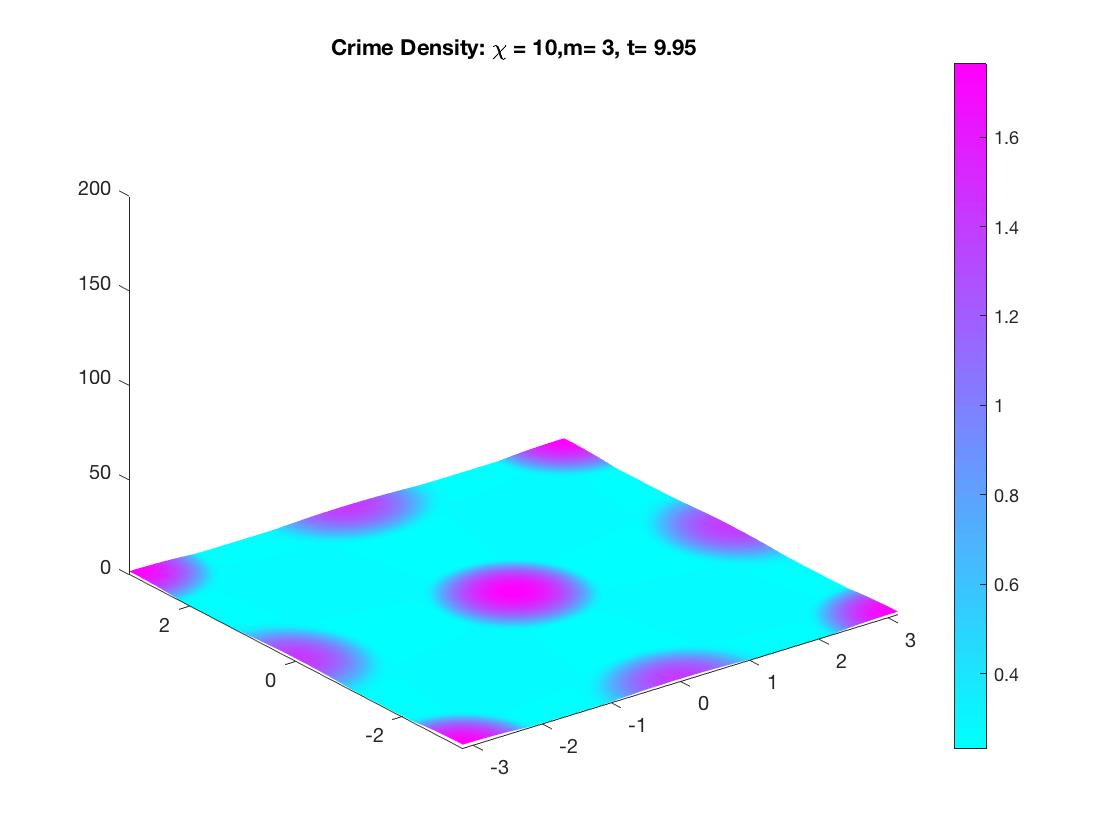}}
   \caption{Numerical solutions with $m=3,\;\chi=10$ and $u =v = \frac{1}{\sqrt{2\pi \sigma^2}}e^{-\frac{|x|^2}{2\sigma^2}}$ with $\sigma = 1/4$. }   \label{fig:oc}
\end{figure}

We are now prepared to construct a solution of (\ref{0}) by means of appropriate compactness arguments.
\begin{lem}\label{lem67}
  Let $m>\frac{3}{2}$.
  Then there exist $(\eps_j)_{j\in\N} \subset (0,1)$ as well as functions
  \be{67.1}
	\left\{ \begin{array}{l}
	u \in L^\infty_{loc}(\bom\times [0,\infty))
	\qquad \mbox{and} \\[1mm]
	v \in C^0(\bom\times [0,\infty)) \cap \bigcap_{q>2} L^\infty_{loc}([0,\infty);W^{1,q}(\Omega))
	\end{array} \right.
  \ee
  such that $\eps_j\searrow 0$ as $j\to\infty$, that $u\ge 0$ a.e.~in $\Omega\times (0,\infty)$ and and $v>0$ in $\bom\times [0,\infty)$,
  that as $\eps=\eps_j\searrow 0$ we have
  \begin{eqnarray}
	& & \ueps\to u
	\qquad \mbox{in } \bigcap_{p\ge 1} L^p_{loc}(\bom\times [0,\infty))
	\mbox{ and a.e.~in } \Omega\times (0,\infty), 
	\label{67.3}
	\\[1mm]
	& & \veps \to v
	\qquad \mbox{in } C^0_{loc}(\bom\times [0,\infty))
	\qquad \mbox{and}
	\label{67.4}
	\\[1mm]
	& & \nabla\veps\wsto \nabla v
	\qquad \mbox{in } \bigcap_{q>2} L^\infty_{loc}([0,\infty);L^q(\Omega)),
	\label{67.5}
  \eea
  and that $(u,v)$ form a global weak solution of (\ref{0}) in the sense of Definition \ref{dw}.
\end{lem}
\proof
  We take any $\alpha>0$ such that
  \bas
	\alpha \ge \left\{ \begin{array}{ll}
	\displaystyle
	\max \Big\{ \frac{m+1}{2} \, , \, 2 \Big\} \equiv 2
	\qquad & \mbox{if } m\in \Big(\frac{3}{2},2\Big], \\[2mm]
	\displaystyle
	m-1
	& \mbox{if } m>2,
	\end{array} \right.
  \eas
  and note that then Lemma \ref{lem51} and Lemma \ref{lem52} may simultaneously be applied so as to show that thanks to Lemma \ref{lem1},
  \bas
	\Big( (\ueps+\eps)^\alpha\Big)_{\eps\in (0,1)}
	\mbox{ is bounded in } L^2((0,T);W^{1,2}(\Omega))
	\quad \mbox{for all } T>0
  \eas
  and that
  \bas
	\Big( \partial_t (\ueps+\eps)^\alpha \Big)_{\eps\in (0,1)}
	\mbox{ is bounded in } L^1\big( (0,T);(W^{2,2}(\Omega))^\star \big)
	\quad \mbox{for all } T>0.
  \eas
  Therefore, employing an Aubin-Lions lemma (\cite{temam}) yields $(\eps_j)_{j\in\N} \subset (0,1)$ and a nonnegative function
  $u$ on $\Omega\times (0,\infty)$ such that $\eps_j\searrow 0$ as $j\to\infty$, and that as $\eps=\eps_j\searrow 0$ we have
  $(\ueps+\eps)^\alpha \to u^\alpha$ in $L^2_{loc}(\bom\times [0,\infty))$ and a.e.~in $\Omega\times (0,\infty)$, whence in particular
  also $\ueps\to u$ a.e.~in $\Omega\times (0,\infty)$.
  Since furthermore Lemma \ref{lem6} warrants boundedness of $(\ueps)_{\eps\in (0,1)}$ in $L^\infty(\Omega\times (0,T))$ for all
  $T>0$, (\ref{67.3}) as well as the inclusion $u\in L^\infty_{loc}(\bom\times [0,\infty))$ result from this due to the Vitali
  convergence theorem.\abs
  As, apart from that, given $T>0$ we know from Lemma \ref{lem666} and Lemma \ref{lem56} that $(\veps)_{\eps\in (0,1)}$ is bounded
  in $C^{\vartheta,\frac{\vartheta}{2}}(\bom\times [0,T])$ and in $L^\infty((0,T);W^{1,q}(\Omega))$ for some 
  $\vartheta=\vartheta(T)\in (0,1)$ and each $q>2$, in view of the Arzel\'a-Ascoli theorem and the Banach-Alaoglu theorem we may
  assume upon passing to a subsequence if necessary that, in fact, $(\eps_j)_{j\in\N}$ is such that with some function $v$
  complying with (\ref{67.1}) we also have (\ref{67.4}) and (\ref{67.5}) as $\eps=\eps_j\searrow 0$.
  The positivity of $v$ in $\bom\times [0,\infty)$ therefore is a consequence of Lemma \ref{lem01}, whereas, finally, the integral
  inequalities in (\ref{w4}) and (\ref{w5}) can be verified in a straightforward manner by relying on (\ref{67.3})-(\ref{67.5})
  when taking $\eps=\eps_j\searrow 0$ in the corresponding weak formulations associated with (\ref{0eps}).
\qed
Our main result on global solvability has thereby actually been established already:\abs
\proofc of Theorem \ref{theo68}. \quad
  All statements have actually been covered by Lemma \ref{lem67} already.
\qed
According to our preparations, and especially due to our efforts to control the dependence of our estimates from
Lemma \ref{lem6} and Lemma \ref{lem56} on $T$ through (K), also the claimed boundedness features can now be obtained
as simple consequences:\abs
\proofc of Theorem \ref{theo69}. \quad
  Again taking the global weak solution of (\ref{0}) obtained in Lemma \ref{lem67}, we only need to observe that thanks to
  the hypotheses (\ref{B1}) and (\ref{B2}), Lemma \ref{lem6} and Lemma \ref{lem56} in conjunction with our notational 
  convention concerning the property (K) guarantee boundedness of $(\ueps)_{\eps\in (0,1)}$ in 
  $L^\infty(\Omega\times (0,\infty))$ and of $(\veps)_{\eps\in (0,1)}$ in $L^\infty((0,\infty);W^{1,q}(\Omega))$ for each
  $q>2$.
  Therefore, namely, the additional features (\ref{69.1}) and (\ref{69.2}) directly result from (\ref{67.3}) and (\ref{67.5}).
\qed
\mysection{Numerical Experiments}\label{sect_numerics}
The purpose of this section is to firstly illustrate how 
the overcrowding effect included in \eqref{0} results in the relaxation of solutions,
and to secondly provide some comparison of this to the situation corresponding to the linear diffusion case $m=1,$ which was not addressed
by our previous analysis.
To this end, we consider the associated evolution problems under initial conditions involving the mildly concentrated data given by
\[
	u_0(x) =v_0(x) = \frac{1}{\sqrt{2\pi \sigma^2}}e^{-\frac{|x|^2}{2\sigma^2}},
	\qquad x\in\Omega,
\]
with some small $\sigma$ on the square $\Omega:=(-3,3)^2$.\abs
We first solve the \eqref{0} numerically with $m=1$ (leading to unconditional diffusion)
and $m=3$ (leading to porous medium type diffusion) with $\sigma = 1/4$.  We illustrate our results for $\chi = 10$ in both simulations, but all other terms are as in the original model proposed in \cite{Short2008} with $B_1 = 1$ and $B_2 = 1$.  
The initial condition for $u$ is illustrated in Figure \ref{fig:init}.  
In the case when $m=1$ we see a concentration of mass around $t= .95$ -- see Figure \ref{fig:blowup}.  
Here there is a real possibility that blow-up happens in finite time, although to make this more precise, more thorough numerical experiments need to be run, which goes beyond the scope of the present work.
What is evident is the concentration around the origin (even if there were 
eventual relaxation) in finite time.  
On the other hand, the power-law diffusion suppresses this concentration entirely as can be observed in Figure \ref{fig:oc}.
The same initial conditions were used
for the numerical experiment illustrated in Figure \ref{fig:oc} with $m=3$.  We clearly see that there is never a concentration of density, and that by time $t=10$, solution comfortably reaches an equilibrium, which may be interpreted as describing
crime hotspots that have spontaneously emerged due to the reaction-cross-diffusion interplay in (\ref{0}). 
Videos of the full simulations can be found in the supplementary material.  These preliminary results lead us to believe that there is 
blow-up when $\chi$ is sufficiently large in the presence of linear diffusion, even for some initial data which are 
only mildly concentrated. However, the considered nonlinear diffusion enhancement suppresses this potential blow-up.

\begin{figure}[H]
  \center
          \subfloat[$m=1,\; t=0$]{\label{fig:1a}\includegraphics[width=0.35\textwidth]{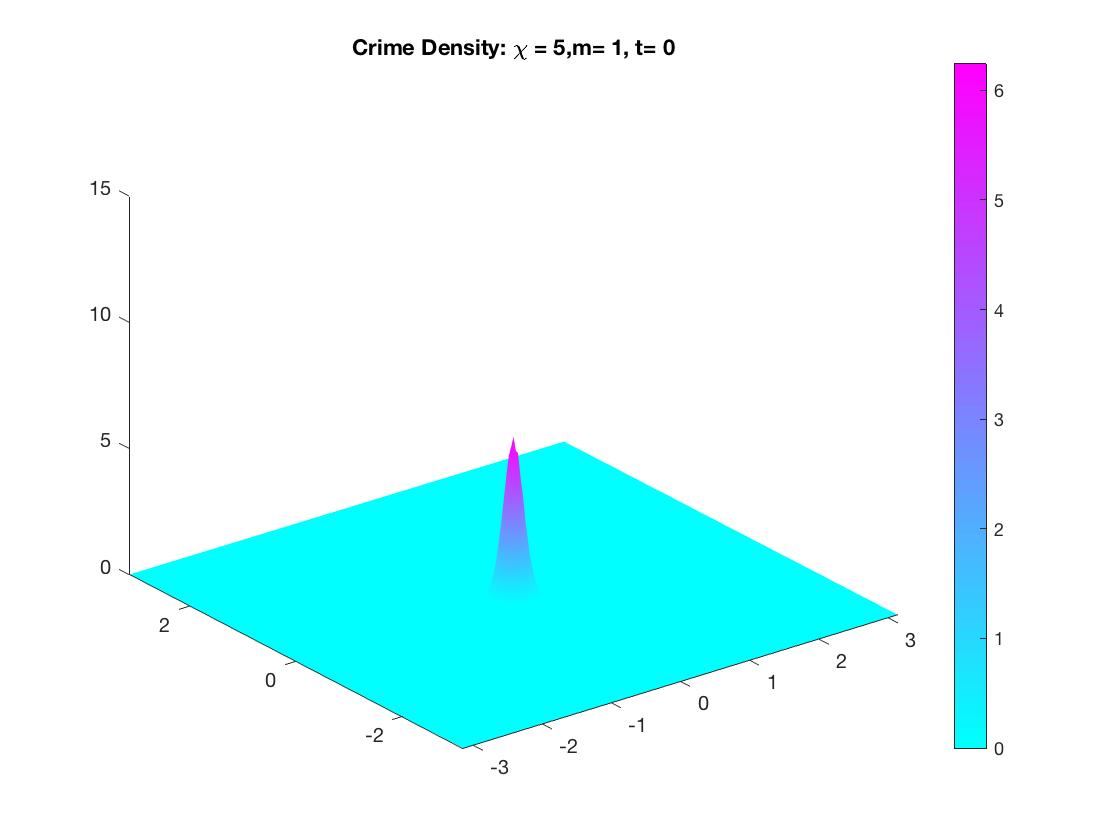}}
     \subfloat[$m=1,\;t = .1$]{\label{fig:2a}\includegraphics[width=0.35\textwidth]{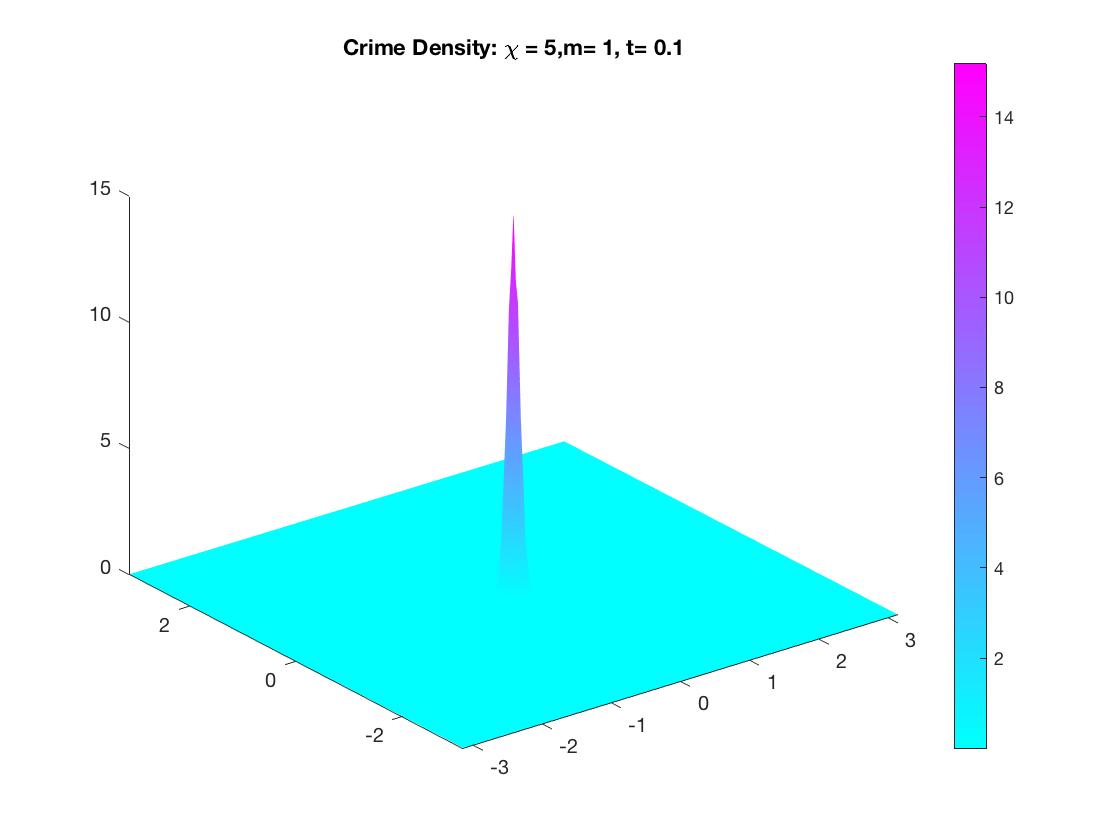}}
  \subfloat[$m=1,\;t=.5$]{\label{fig:3a}\includegraphics[width=0.35\textwidth]{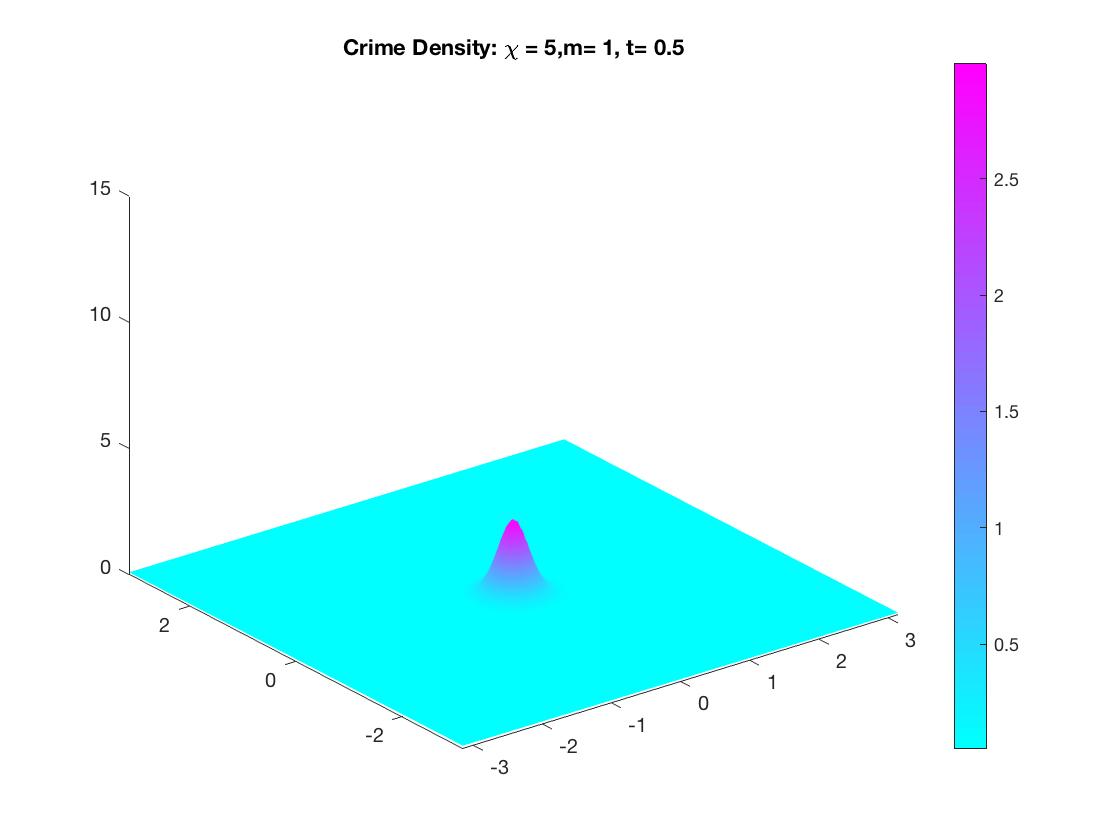}}\\
            \subfloat[$m=3,\; t=0$]{\label{fig:1b}\includegraphics[width=0.35\textwidth]{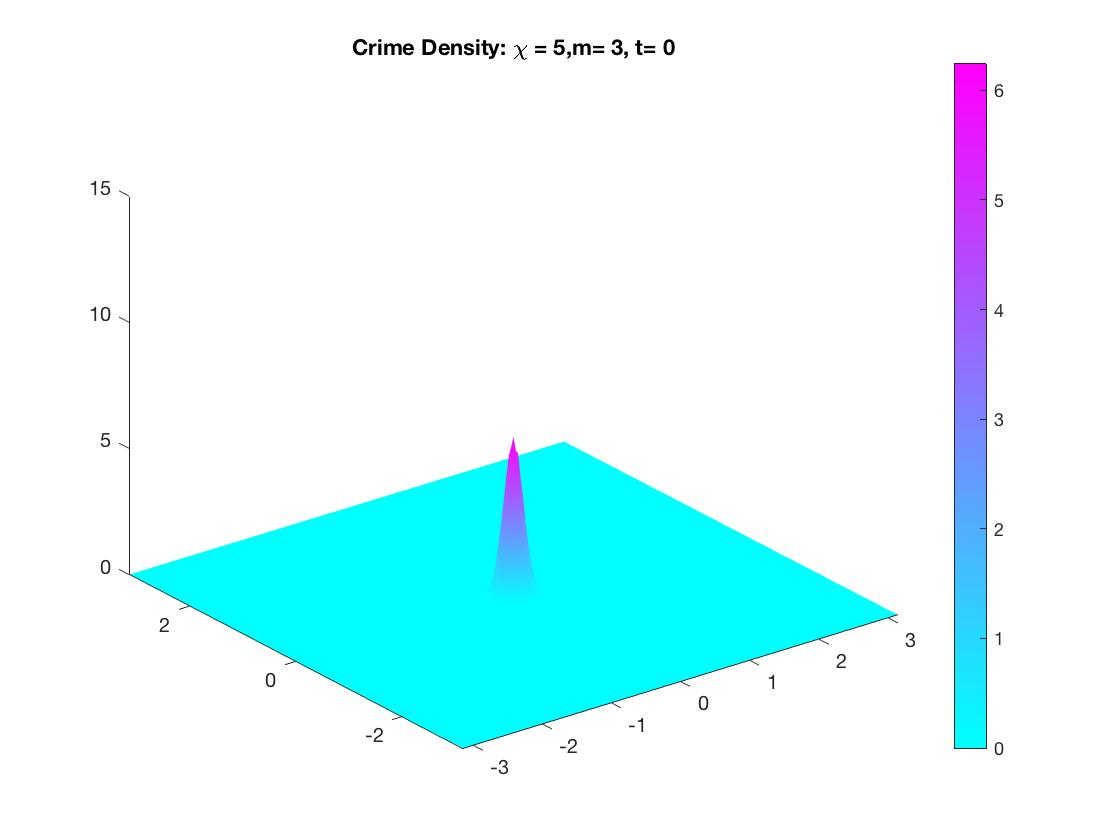}}
     \subfloat[$m=3,\;t = .1$]{\label{fig:2b}\includegraphics[width=0.35\textwidth]{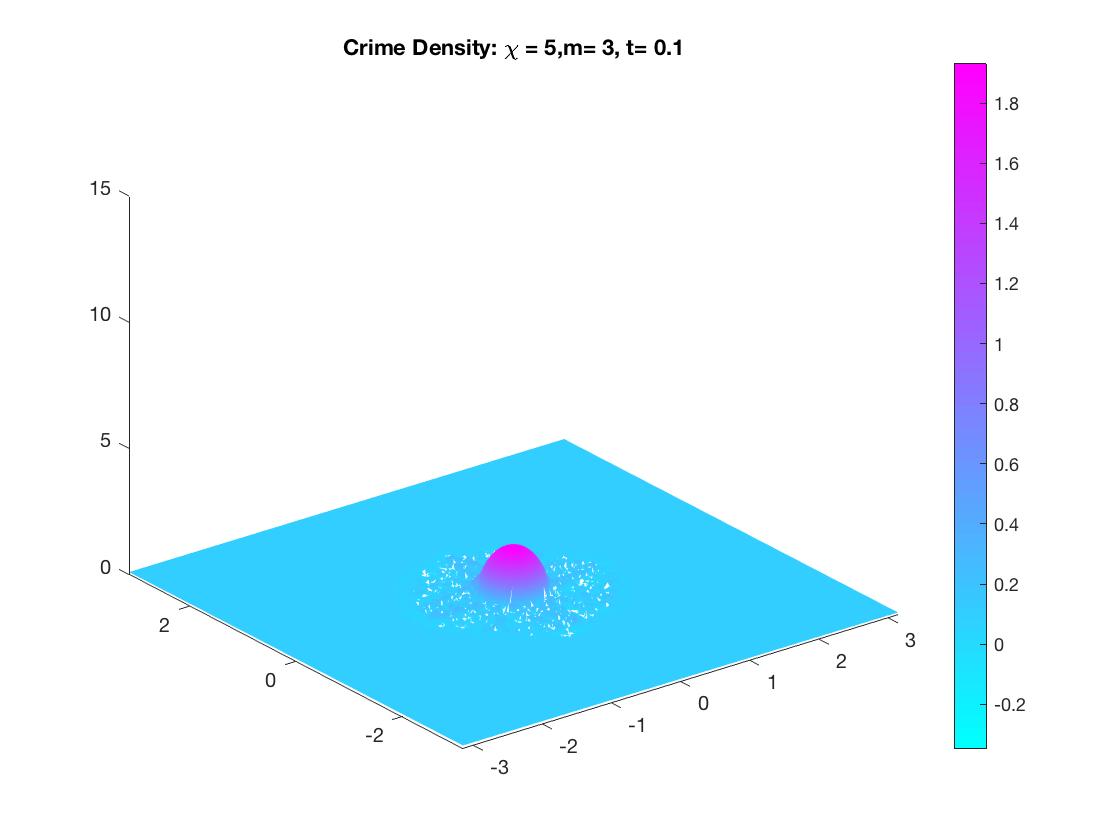}}
  \subfloat[$m=3,\;t=.5$]{\label{fig:3b}\includegraphics[width=0.35\textwidth]{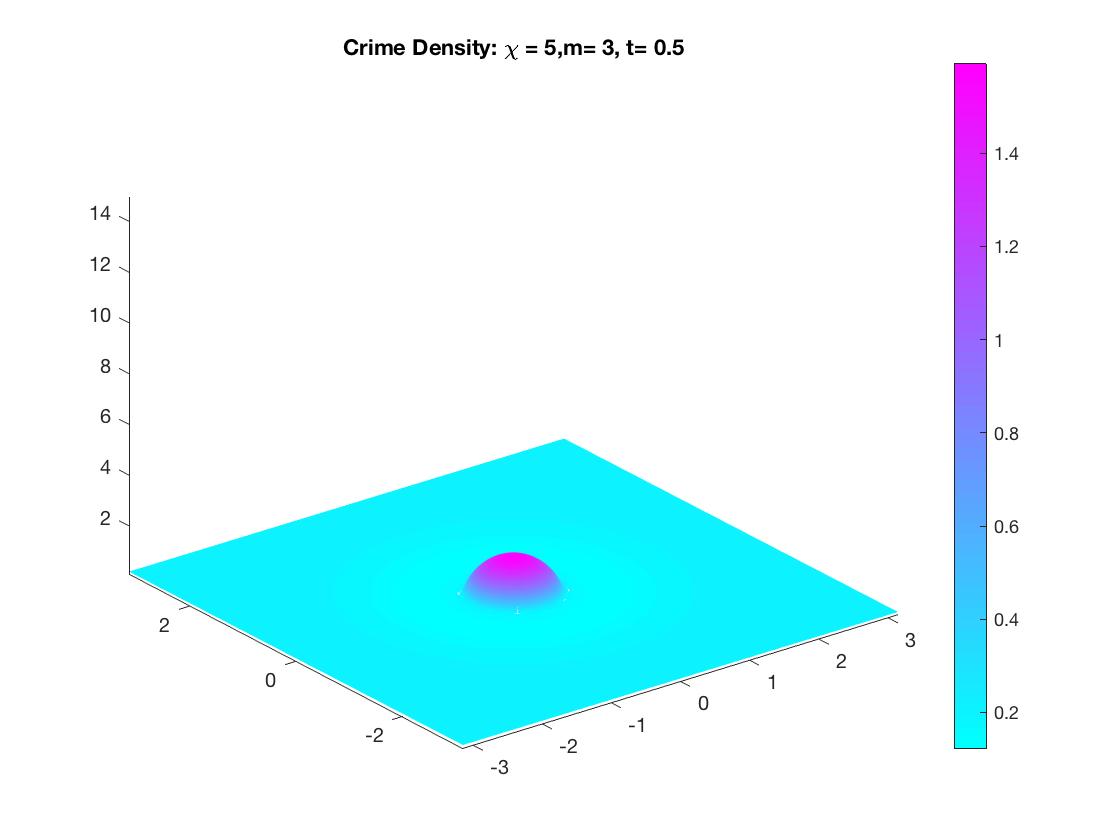}}\\
  \caption{Numerical solutions comparing $m=1$ and $m=3$ with $\chi=5$ and $u =v = \frac{1}{\sqrt{2\pi \sigma^2}}e^{-\frac{|x|^2}{2\sigma^2}}$ with $\sigma = .16$. }   \label{fig:oc1}
\end{figure}

From more general numerical experiments, we observe that the smaller $\chi$ is, the more concentrated the initial data needs to be in order for a potential blow-up to occur. 
Moreover, in the case $m=1$, for each $\chi$ there are initial data which are not concentrated enough to see blow-up, but concentrated enough to see initial growth.  
This initial growth is suppressed by the overcrowding effect seen in \eqref{0}.  This is shown in Figure \ref{fig:oc1} where the top row illustrates the linear diffusion case ($m=1$)
and the bottom row illustrates the non-linear diffusion with $m=3$.  In the top row we observe the initial growth of the solution in Figure \ref{fig:2a}.  This growth does not last for very long
and the solution is already decaying at time $t=.5$ as illustrated in Figure \ref{fig:3a}.  We do see some numerical instabilities for the case $m=3$ around the region of concentration.  We expect that this
is due to the degeneracy of the diffusion and more sophisticated numerical methods need to be used to deal with potential contact lines. 

\bigskip

{\bf Acknowledgement.} \quad
The first author acknowledges support the National Science Foundation, NSF DMS-1516778.  The second author acknowledges support of the {\em Deutsche Forschungsgemeinschaft} in the context of the project 
  {\em Emergence of structures and advantages in cross-diffusion systems} (No.~411007140, GZ: WI 3707/5-1).

\bigskip

\end{document}